\documentclass[11pt]{amsart}

\numberwithin{equation}{section} 
\usepackage{amsmath,amscd}
\usepackage{bm,epsfig,ifthen}
\usepackage{amsfonts,amssymb}
\usepackage{colortbl}
\usepackage{pstricks}
\usepackage{mathtools}
\usepackage{algorithm,algorithmic}

\usepackage[font=small,labelfont=bf]{caption}

\usepackage{color}
\usepackage{graphicx,graphics}
\usepackage{amsfonts}

\usepackage{sidecap}
\graphicspath{{./figures/},{./Figures/},{./Figs_estimation/}}

\usepackage{booktabs}

\usepackage[colorlinks=true, pdfstartview=FitV, linkcolor=blue,citecolor=blue, urlcolor=blue]{hyperref}  

\DeclarePairedDelimiterX\braket[2]{\langle}{\rangle}{#1 \delimsize\vert #2}

\newcommand{\norm}[1]{\left\lVert#1\right\rVert}
\usepackage{newtxtext}

\numberwithin{equation}{section} \oddsidemargin=-.0cm
\newtheorem{thm}{Theorem}[section]
\newtheorem{prop}{Proposition}[section]
\newtheorem{lem}{Lemma}[section]
\newtheorem{defi}{Definition}[section]

\newtheorem{thm_app}{Theorem}
\newtheorem{prop_app}[thm_app]{Proposition}

\newtheorem{rem}{Remark}[section]
\newtheorem{rem_app}{Remark}
\newtheorem{cor}{Corollary}[section]

\def\bt{\begin{thm}}
\def\bprop{\begin{prop}}
\def\eprop{\end{prop}}

\def\et{\end{thm}}
\def\bl{\begin{lem}}
\def\el{\end{lem}}
\def\bd{\begin{defi}}
\def\ed{\end{defi}}
\def\bc{\begin{cor}}
\def\ec{\end{cor}}
\def\bp{\begin{proof}}
\def\ep{\end{proof}}
\def\br{\begin{rem}}
\def\er{\end{rem}}

\def\Forall{\text{ } \forall \:}
\def\d{\, \mathrm{d}}

\def\Res{\, \mathrm{Re}\,}
\def\Ims{\, \mathrm{Im}\,}
\def\T{\, \widehat{T}}

\def\bi{\begin{itemize}}
\def\ei{\end{itemize}}
\def\be{\begin{equation}}
\def\ee{\end{equation}}
\def\bes{\begin{equation*}}
\def\ees{\end{equation*}}
\def\bea{\begin{equation} \begin{aligned}}
\def\eea{\end{aligned} \end{equation}}
\def\beas{\begin{equation*} \begin{aligned}}
\def\eeas{\end{aligned} \end{equation*}}

\def\Xi{\X}
\def\s{\mathfrak{s}}
\def \c{\mathfrak{c}}

\def\y{\boldsymbol{y}}
\def\x{\boldsymbol{x}}
\def\z{\boldsymbol{z}}
\def\w{\boldsymbol{w}}

\def\X{\boldsymbol{X}}
\def\cH{\mathcal H}

\begin{document}

\title[Effective Reduced Models from DDEs: Bifurcations, Tipping Solution Paths, and ENSO variability]{Effective Reduced Models from Delay Differential Equations:\\
Bifurcations, Tipping Solution Paths, and ENSO variability}
\author{Micka{\"e}l D. Chekroun}
\address[MC]{Department of Earth and Planetary Sciences, Weizmann Institute, Rehovot 76100, Israel\\
Department of Atmospheric and Oceanic Sciences, University of California, Los Angeles, CA 90095-1565, USA}
\email{mchekroun@atmos.ucla.edu}

\author{Honghu Liu}
\address[HL]{Department of Mathematics, Virginia Tech, Blacksburg, VA 24061, USA} 
\email{hhliu@vt.edu}

\date{January 2, 2024}

\begin{abstract}
Conceptual delay models have played a key role in the analysis and understanding of El Ni\~no-Southern Oscillation (ENSO) variability. 
Based on such delay models,  we propose in this work a novel scenario  for the fabric of ENSO variability resulting from the subtle  interplay between stochastic disturbances and nonlinear invariant sets emerging from bifurcations of the unperturbed dynamics. 

To identify these invariant sets we adopt an approach combining Galerkin-Koornwinder (GK) approximations of delay differential equations and center-unstable manifold reduction techniques.   In that respect, GK approximation formulas are  
reviewed and synthesized, as well as analytic approximation formulas of center-unstable manifolds.
The reduced systems derived thereof enable us  to conduct a thorough analysis of the bifurcations arising in a standard delay model of ENSO.   We identify thereby a saddle-node bifurcation of periodic orbits co-existing with a subcritical Hopf bifurcation, and a homoclinic bifurcation for this model. We show furthermore that the computation of unstable periodic orbits (UPOs) unfolding through these bifurcations is considerably simplified from the reduced systems.

These dynamical insights enable us in turn to design a stochastic model whose solutions---as the delay parameter drifts slowly through its critical values---produce a wealth of temporal patterns resembling ENSO events and exhibiting also decadal variability. 
Our analysis dissects the origin of this variability and shows how it is tied to certain transition paths between invariant sets of the unperturbed dynamics (for ENSO's interannual variability) or simply due to the presence of UPOs close to the homoclinic orbit (for decadal variability).  In short, this study points out the role of solution paths evolving through tipping ``points'' beyond equilibria, as possible mechanisms organizing  the variability of certain climate phenomena. 


\end{abstract}

\maketitle

\tableofcontents

\section{Introduction}
Since its inception in the 80s and 90s as conceptual models to study the El Ni{\~n}o-Southern Oscillation (ENSO) phenomenon  \cite{Suarez_al88,BH89,TSCJ94,Neelin_al98}, delay models have attracted a growing attention in climate modeling; see  \cite{Galanti_al00,ghil2008delay,koren2011aerosol,roques2014parameter,GCStep15,koren2017exploring,keane2017climate,CGN18,boers2017inverse,falkena2019derivation}. 
However, only recently the bifurcation analysis of these conceptual delay models arising in climate has been undertaken with modern tools dedicated to delay models; see e.g.~\cite{krauskopf2014bifurcation,keane2015delayed,keane2017,keane2019effect,CKL20}.

 In this work, we pursue this endeavor by providing a detailed analysis of the bifurcations arising in the ENSO delay model of Suarez and Schopf \cite{Suarez_al88}. Partial bifurcation results have been reported about this model. These results were often treated in a subsidiary fashion within a work of more general scope 
\cite{falkena2019derivation,anikushin2023hidden} relegating details about certain calculations and dynamical insights to the reader. Here, we provide a thorough analysis of these bifurcations and associated calculations to help the reader build up intuitions in view of applications, in particular regarding the stochastic modeling of ENSO as discussed also in this study.  
Our approach employs the Galerkin method based on the Koornwinder basis functions \cite{Koo84} as introduced in  \cite{CGLW16} to derive rigorous low-dimensional ordinary differential equation (ODE) approximations (GK systems for short), to which a reduction to the center-unstable manifold is applied, following \cite{CKL20}. Here, both the derivation of GK systems and the center-unstable manifold calculations are made self-contained and more transparent than in  \cite{CGLW16,CKL20} to better serve the ``practitioner''  interested in applications.  

For the Suarez and Schopf delay model, the resulting two-dimensional reduced ODE system demonstrates remarkable skills in predicting the local and global bifurcations including a subcritical Hopf bifurcation, a homoclinic bifurcation, and a saddle-node bifurcation of periodic orbits (SNO bifurcation) when the delay parameter $\tau$ is varied.  Following \cite[Theorem III.1]{CKL20}, the subcritical Hopf bifurcation is proved by analyzing the sign of the Lyapunov coefficient  
of the reduced system.  Theorem III.1 of \cite{CKL20} simplifies the calculation of this coefficient to that of basic operations consisting of solving  
linear algebraic systems with triangular matrices, computing the GK spectrum at the critical delay parameter, and forming the involved inner products;   see \eqref{Eq_l1_GK} below.
The detection of homoclinic and SNO bifurcations benefits from the stable and accurate approximation skills of the Unstable Periodic Orbits (UPOs) to the Suarez and Schopf model, by simple backward integration of the reduced system. The computation of UPOs is known to be challenging for high-dimensional problems and to require a careful formulation and implementation of Newton methods and the like \cite{gritsun2008unstable,gritsun2013statistical}.  Here, our reduced system allows thus for bypassing these difficulties.  

 There is a relatively vast body of literature dealing with Galerkin methods exploiting other basis functions to approximate delay model by ODEs, including step functions \cite{Banks_al78,kappel1978autonomous}, splines \cite{banks1979spline,Banks_al84}, linear and sine functions \cite{Vyasarayani12,Wahi_al05}, and Legendre polynomials \cite{Kappel86,Ito_Teglas86}. Techniques based on pseudospectral collocation methods have also emerged as an alternative route to approximate delay differential equations (DDEs) by ODEs \cite{breda2016pseudospectral} and have shown relevance in computing characteristic roots  \cite{breda2005pseudospectral}, Lyapunov exponents \cite{breda2014approximating}, and bifurcation analysis \cite{ando202215} using  ODE bifurcation software packages such as AUTO \cite{Auto} and MatCont \cite{dhooge2008new,Matcont}.  
These methods, although showing great promises, do not have yet developed into toolboxes such as  DDE-BIFTOOL \cite{Engelborghs_al02,sieber7144dde} and KNUT \cite{Knut} which provide currently the most evolved software packages to analyze bifurcations arising in DDEs with discrete delays.
Nevertheless, ODE approximation methods offer a framework that is more flexible to analyze bifurcations, covering a broader range of applications as encompassing delay models with (possibly nonlinear) functional of distributed delays \cite{CGLW16},  or renewal equations \cite{breda2016pseudospectral}.  

Yet, the effective derivation of reduced systems from the ODE approximations produced by such methods, seems to have been poorly exploited.
As shown in the case of GK systems, to dispose of such reduced systems has many potential attributes that await to be 
tested in applications involving DDEs. Steps in this direction include e.g.~the effective computation of low-dimensional solutions to (nearly) optimal controls in feedback form, enabling to avoid solving the infinite-dimensional Hamilton-Jacobi-Bellman (HJB) equation associated with the optimal control problem of the original DDE \cite{CKL18_DDE}. 

Reduced models derived from GK systems have also shown recently their relevance to improve the realism of DDE's  solutions via stochastic modeling, as illustrated in cloud physics \cite{Chekroun_al22SciAdv}. 
Conceptual delay models were proposed in \cite{koren2011aerosol,koren2017exploring} to envision open cellular cells formed by marine stratocumulus clouds as an oscillatory predator-prey mechanism of clouds (prey)
scavenged by rain (predator). Such conceptual models produce oscillations that, albeit grounded on physical intuitions,  are too regular (periodic) to bear the comparison with e.g.~satellite observations or high-end model simulations.
In \cite{Chekroun_al22SciAdv}, a significant step was taken  to bridge the gap between the latter and these conceptual models. To do so, the reduced systems derived in \cite{Chekroun_al22SciAdv} from GK approximations of the DDE model of  \cite{koren2011aerosol} allowed for identifying 
nonlinear structures providing the level set of the oscillation's phase function, known as isochrons \cite{guckenheimer1975isochrons,Ashwin2016}.
Current algorithms to compute isochrons  \cite{mauroy2014global,detrixhe2016fast} suffer the ``curse of dimensionality,'' making the direct computation of isochrons out of reach, for DDEs. 
 In  \cite{Chekroun_al22SciAdv}, the knowledge of the isochrons in the reduced state space was shown to be sufficient to  
design stochastic parameterizations interacting with the genuine DDE's isochrons, leading in turn to  stochastic chaotic solutions with new virtues including enhanced time-variability mimicking that of clouds' oscillations.

In this study, benefiting from the high-accuracy approximation skills achieved by our 2D reduced GK system, we reach a detailed understanding of the DDE phase portrait when the delay parameter  $\tau$ is varied; see Fig.~\ref{Fig_bif_combo} below. This understanding allows us to dissect the response of the Suarez and Schopf model when subject to stochastic disturbances while the delay parameter drifts slowly through the aforementioned bifurcations.  As shown in Sec.~\ref{Sec_nln_blocks}, we obtain indeed tipping solution path (TSP) whose certain time episodes resemble ENSO-like temporal patterns that are explained in terms of transition paths between invariant sets of the DDE's unperturbed dynamics, such as locally stable/unstable steady states and UPOs.  By conducting in Sec.~\ref{Sec_homocline} a spectral analysis of the TSP's frequency content,  we show furthermore that the TSP's irregular behavior on longer timescales is dominated by decadal variability such as documented in recent ENSO studies \cite{dieppois2021enso}. The dynamical insights gained in Sec.~\ref{Sec_Reduc_ENSO} allow us to trace back the origin of this decadal variability in our stochastic model which results from the presence of UPOs located close to the homoclinic orbit.

The remainder of this paper is organized as follows. We first recall the key analytic and algebraic elements for the formal  construction of the Galerkin-Koornwinder (GK) approximations of general scalar DDEs in Sec.~\ref{GK_section}. We survey in Sec.~\ref{Sec_GK_centerman} analytic approximation formulas of center-unstable manifolds---including leading-order  (Theorem~\ref{Thm_CM_approx}) and higher-order formulas (Sec.~\ref{Sec_Higher-order})---that we present for the reduction of GK systems experiencing a loss of stability at a critical value of the delay parameter $\tau$. 
These formulas are then applied to the Suarez and Schopf model in Sec.~\ref{Sec_Reduc_ENSO} to derive an effective 2D reduced GK system with $\tau$-dependent coefficients (Eq.~\eqref{Eq_EffectiveReduced_ENSO}). 
Sections \ref{Sec_Bif_Hopf} and \ref{Sec_Bif_SNO} present the reduced system skills in predicting for the delay model a subcritical Hopf bifurcation, and an SNO bifurcation co-existing with an homoclinic bifurcation, respectively. 
Section \ref{Sec_approx_orbits} provides numerical evidences of the  highly-accurate approximation skills achieved by the reduced systems, in particular regarding the computation of UPOs and limit cycles. Section \ref{Sec_model_error} gives then model error estimates complementing these numerical results. 
Finally, in Sec.~\ref{Sec_ENSO_var}, we take  advantage  of the dynamical insights gained in Sec.~\ref{Sec_Reduc_ENSO} to design a stochastic model whose solutions exhibit ENSO-like patterns and decadal variability. Some    final remarks and potential future directions are then outlined in Sec.~\ref{Sec_conclusion}.

\section{Galerkin-Koornwinder (GK) approximations of DDEs}\label{GK_section}

We consider nonlinear {\it scalar} DDEs of the form
\be \label{Eq_DDE}
\frac{\d x(t)}{\d t} = a x(t) + b x(t-\tau) + c \int_{t-\tau}^t x(s)\d s + F\Big(x(t), x(t-\tau),\int_{t-\tau}^t x(s) \d s\Big),
\ee
where $a$, $b$ and $c$ are real numbers, $\tau> 0$ is the delay parameter, and $F$
is a nonlinear function. 
We restrict ourselves to the scalar case to simplify the  presentation, but the approach extends 
to systems of nonlinear DDEs involving possibly several delays as detailed in \cite{CGLW16} and illustrated in \cite{CKL20} for the cloud-rain delay model of \cite{koren2011aerosol}.

\subsection{Koornwinder polynomials}
First, let us recall that Koornwinder polynomials $K_n$ are obtained from 
Legendre polynomials $L_n$, for any nonnegative integer $n$, according to the relation
\be \label{eq:Pn}
K_n(s)= -(1+s)\frac{\d}{\d s} L_n(s) +( n^2 + n + 1) L_n(s), \quad s \in [-1, 1],
\ee
see \cite[Eq.~(3.3)]{CGLW16}.

Koornwinder polynomials are  known to form an orthogonal set for the following weighted inner product on $[-1,1]$ with a point-mass, $\mu(\d x)= \frac{1}{2} \d x  + \delta_{1},$
where $\delta_1$ denotes the Dirac point-mass at the right endpoint $x=1$; see  \cite{Koo84}.  In other words, the following orthogonality property holds:
\bea
\int_{-1}^{1} K_n(s) K_p(s) \d \mu (s)& = \frac{1}{2} \int_{-1}^{1} K_n(s) K_p(s) \d s + K_n(1) K_p(1)\\
&=0, \, \mbox{ if } p\neq n.
\eea

It is also worthwhile noting that Koornwinder polynomials augmented by the right endpoint values as follows 
\be \label{eq:Pn_prod}
\mathcal{K}_n = (K_n, K_n(1)),
\ee
form an orthogonal basis of the  product space 
\bes \label{eq:E}
\mathcal{E} = L^2([-1,1); \mathbb{R}) \times  \mathbb{R},
\ees
endowed with the inner product:
\be \label{eq:inner_E}
\Big\langle (f, p), (g, q) \Big\rangle_{\mathcal{E}}  = \frac{1}{2} \int_{-1}^1 f(s)g(s) \d s  + pq, \qquad \text{for all }  (f, p), (g, q) \in \mathcal{E}.
\ee
The norm induced by this inner product is denoted by $\|\cdot\|_{\mathcal{E}}$. The basis function $\mathcal{K}_n$ has then its $\|\cdot\|_{\mathcal{E}}$-norm  given by \cite[Prop.~3.1]{CGLW16}:
\be \label{eq:Kn_norm}
\|\mathcal{K}_n\|_{\mathcal{E}} = \sqrt{\frac{(n^2+1)((n+1)^2+1)}{2n+1}}.
\ee
This is a useful property for calculating the GK approximations, as it will come apparent below.

Finally, since the original Koornwinder basis is given on the interval $[-1, 1]$ and the state space of a DDE such as Eq.~\eqref{Eq_DDE} is made of functions defined on $(-\tau,0)$, we will work mainly with the following rescaled version of  Koornwinder polynomials. 
The rescaled Koornwinder polynomials $K_n^\tau$ are defined as follows 
\bea \label{eq:Pn_tilde}
K^\tau_n\colon  [-\tau, 0] & \rightarrow \mathbb{R}, \\
\theta & \mapsto  K_n \Bigl( 1 + \frac{2 \theta }{\tau} \Bigr).
\eea
They form orthogonal polynomials on the interval $[-\tau, 0]$ for the $L^2$-inner product on $(-\tau,0)$ with a Dirac point-mass at the right endpoint $0$.

The following family of rescaled Koornwinder polynomials augmented with a right endpoint value,
\be \label{eq:Pn_tilde_prod}
\mathcal{K}_n^\tau = (K_n^\tau, K_n^\tau(0)),
\ee
forms then an orthogonal basis for the Hilbert space $\mathcal{H} = L^2([-\tau,0); \mathbb{R}) \times  \mathbb{R}$ endowed with the inner product $\langle \cdot, \cdot \rangle_{\mathcal{H}}$ given by \eqref{eq:inner_E} in which the integral is taken on $(-\tau,0)$ and weighted by $1/\tau$. Note that since $K_n(1)=1$ \cite[Prop.~3.1]{CGLW16}, we have also 
\be\label{Eq_normalization}
K_n^\tau(0)  = 1.
\ee
Finally, observe that $\|\mathcal{K}_n^\tau\|_{\cH} = \|\mathcal{K}_n\|_{\mathcal{E}}$.

\subsection{Space-time representation and GK approximations}
It is well-known that a DDE such as Eq.~\eqref{Eq_DDE} is an infinite-dimensional dynamical system \cite{Hale88} in which a history segment has to be specified over the interval $[-\tau,0)$ to apprehend properly the existence and computation problems of its solutions \cite{Hale_Lunel93,curtain1995,bellenbook}. 
  
\begin{figure}[hbtp]
   \centering
\includegraphics[width=0.6\textwidth, height=0.5\textwidth]{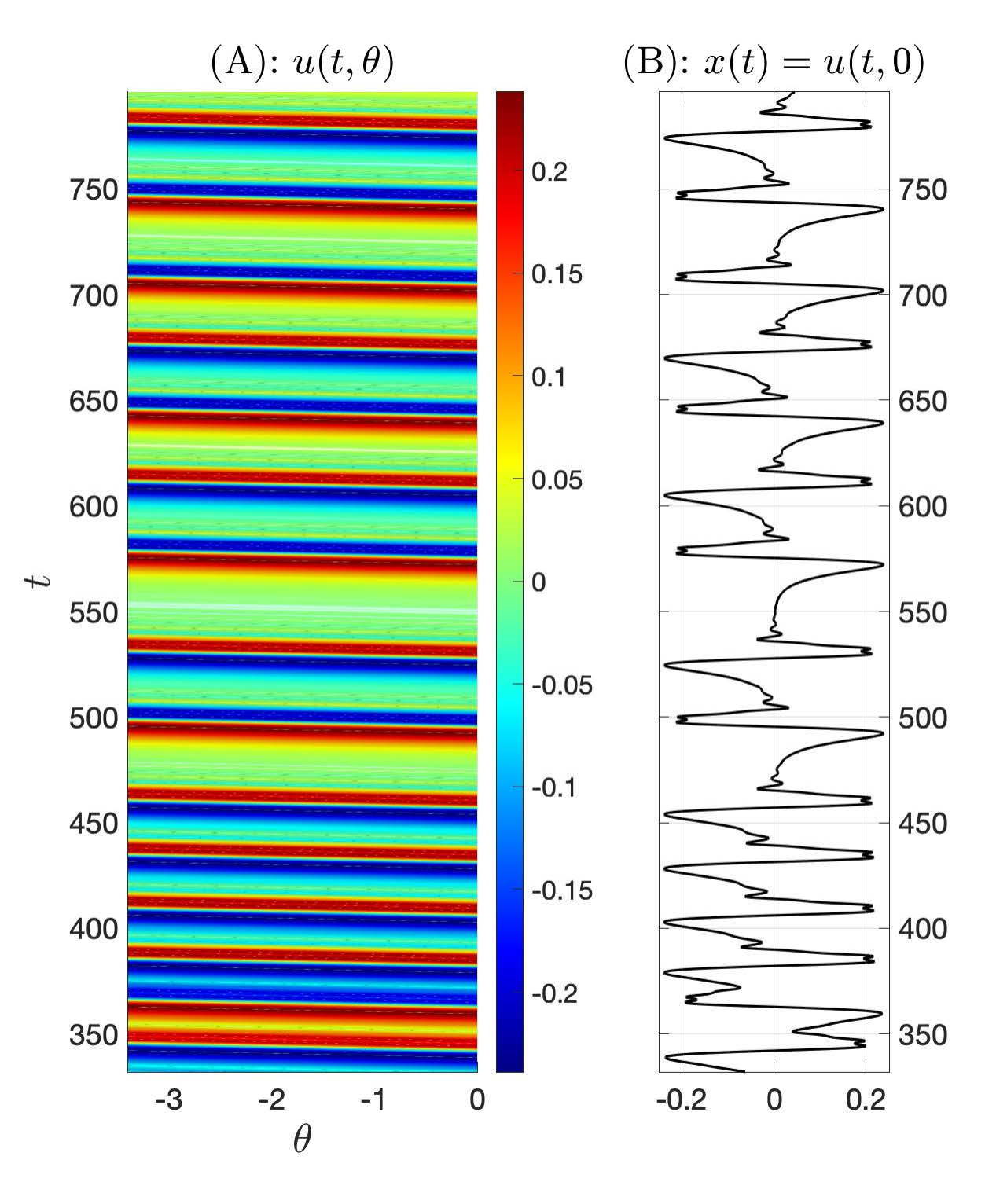}
  \caption{Panel A shows a solution $u(t, \theta)$ to a transport problem  of the form \eqref{lin_PDE}--\eqref{PDE_BC} obtained as a reformulation of the DDE $\dot{x}=a x(t-\tau)-b x(t-\tau)^3$ from \cite{CGLW16}. Panel B shows the time series $u(t,0)$ at the right endpoint $\theta=0$. It coincides with the time series $x(t)$ obtained by solving the original DDE. The solution is shown for $a=0.5$, $b=20$ and $\tau=0.4$.} \label{Fig_hovmoller}
\end{figure}

Thus, given a function of time, $x$, solving Eq.~\eqref{Eq_DDE} one distinguishes between the {\it history segment} $\{x(t + \theta): \theta \in [-\tau, 0)\}$ and the {\it current state}, $x(t)$. Denoting by $u(t,\theta)$ the history segment, one can then rewrite  the DDE  \eqref{Eq_DDE} as the transport equation 
\be\label{lin_PDE}
\partial_t u = \partial_{\theta}u, \quad -\tau \le \theta < 0,
\ee
subject to the {\it nonlocal and nonlinear boundary condition}
\bea\label{PDE_BC}
\partial_{\theta} u|_{\theta = 0}& = a u(t,0) + b u(t,-\tau) + c \int_{t-\tau}^t u(s,0)\d s\\
&\qquad \;\; + F\Big(u(t,0), u(t,-\tau),\int_{t-\tau}^t  u(s,0) \d s\Big), \qquad \mbox{for $t \geq 0$.}
\eea
This reformulation is often called the space-time representation in the literature; see Fig.~\ref{Fig_hovmoller} for an illustration.

Using the rescaled  Koornwinder polynomials \eqref{eq:Pn_tilde}, one can then seek for  approximations of $u(t,\theta)$ solving  Eqns.~\eqref{lin_PDE}-\eqref{PDE_BC} under the form
\be\label{Eq_App1}
u_N(t,\theta) = \sum_{j=0}^{N-1} y_j(t) K_j^\tau (\theta).
\ee
Given the property \eqref{Eq_normalization}, the current state $x(t)=u(t,0)$ is approximated 
by
\be\label{Eq_App2}
x_N(t)= \sum_{j=0}^{N-1} y_j(t). 
\ee
Note that in \eqref{Eq_App1} the Koornwinder polynomials are indexed according to their degree $j$.  

The question arises then of determining the equations that the coefficients $y_j(t)$ in \eqref{Eq_App1} must satisfy in order to have that  $u_N(t,\theta)$ provides a rigorous approximation that converges to $u(t,\theta)$ as $N\rightarrow \infty.$
This nontrivial problem has been solved in \cite{CGLW16}. The next section recalls the structure of these equations forming what we call a GK system.

\vspace{.5cm}
\subsection{The formulas of GK approximations}\label{Sec_GKaction}
The analysis of  \cite{CGLW16} shows that $y_j(t)$ solving the $N$-dimensional system of ODEs \eqref{Galerkin_AnalForm} provides a rigorous method for approximating the solutions to a broad class of DDEs; see \cite[Sec.~5]{CGLW16}. 
Such systems are given by: 
\begin{equation} \label{Galerkin_AnalForm}
\begin{aligned}
\frac{\d y_j}{\d t} & = \frac{1}{\|\mathcal{K}_j\|_{\mathcal{E}}^2 } \sum_{n=0}^{N-1} \bigg( a + b K_n(-1) + c \tau (2 \delta_{n,0} - 1) \\
& \hspace{8em} + \frac{2}{\tau}\sum_{k=0}^{n-1} a_{n,k} \left( \delta_{j,k} \|\mathcal{K}_j\|^2_{\mathcal{E}} - 1 \right) \bigg) y_n(t) \\
& \hspace{1em} + \frac{1}{\|\mathcal{K}_j\|_{\mathcal{E}}^2} F \left( \sum_{n=0}^{N-1} y_n(t),  \sum_{n=0}^{N-1} y_n(t) K_n(-1), \tau y_0(t) - \tau \sum_{n=1}^{N-1} y_n(t) \right),\\
& \;\; 0\leq j\leq N-1.
\end{aligned}
\end{equation}
Here the Kronecker symbol $\delta_{j,k}$ has been used, and the coefficients $a_{n,k}$ are obtained by solving a triangular linear system in which the right-hand side has 
explicit coefficients depending on $n$ \cite[Prop.~5.1]{CGLW16}; see Appendix~\ref{sect:coef_matrix_proof}. 
The values of the $K_n(-1)$'s are known and recalled in \eqref{Eq_K-1} below, and we recall that  the formula for the $\|\mathcal{K}_n\|_{\mathcal{E}}$ is given by \eqref{eq:Kn_norm}.


We rewrite now the GK system \eqref{Galerkin_AnalForm} into the following compact form:
\be \label{Galerkin_cptForm}
\boxed{\frac{\d \y}{\d t} = A_N (\tau) \y + G_2 (\y,\tau),}
\ee
where $\y= (y_0, \cdots, y_{N-1})^{T}$, 
and in which the terms $A_N (\tau) \y$ and $G_2(\y,\tau)$ group the linear and nonlinear terms in Eq.~\eqref{Galerkin_AnalForm}, respectively. 
The entries of the $N\times N$ matrix  $A_N (\tau)$ are given by \cite[Eq.~(5.20)]{CGLW16}
\begin{equation} \label{eq:A}
\boxed{
\begin{aligned}
\Big(A_N (\tau)\Big)_{i,j}& = \frac{1}{\|\mathcal{K}_i\|_{\mathcal{E}}^2 }\Big(a + b K_j(-1) + c \tau (2 \delta_{j,0} - 1) \\
& \hspace{8em}+ \frac{2}{\tau}\sum_{k=0}^{j-1} a_{j,k} \left( \delta_{i,k} \|\mathcal{K}_i\|^2_{\mathcal{E}} - 1 \right ) \Big), \;\; i, j = 0, \cdots, N-1.
\end{aligned}
}
\end{equation}
 Here, the coefficients $a_{j,k}$ are independent of the DDE model and are determined by \cite[Prop.~5.1]{CGLW16}; see Proposition \ref{prop:dPn} in Appendix~\ref{sect:coef_matrix_proof}. 
Re-writing $A_N(\tau)$ as $A_N(\tau)=2/\tau \mathcal{P}_N +\mathcal{Q}_N$, one observes that only $\mathcal{Q}_N$ depends on the model's  parameters. More generally, the matrix $A_N(\tau)$  accounts for the approximation of the transport equation $\partial_t u=\partial_\theta u$ and the contribution of the linear terms arising in the nonlocal boundary condition \eqref{PDE_BC}; see \cite{CGLW16}.

The nonlinear mapping $G_2$ has its components $F_{N}^j$ given by \cite[Eq.~(5.21)]{CGLW16} 
\be \label{eq:G}
\boxed{F_{N}^j(\y,\tau) = \frac{1}{\|\mathcal{K}_{j}\|_{\mathcal{E}}^2} F \left( \sum_{n=0}^{N-1} y_n(t),  \sum_{n=0}^{N-1} y_n(t) K_n(-1),\tau y_0(t) - \tau \sum_{n=1}^{N-1} y_n(t) \right),}
\ee
for each  $0\leq j\leq N-1$.

Due to rigorous convergence results of GK approximations \cite{CGLW16}, the GK formulas above provide a powerful apparatus to analyze DDEs by means of ODE approximations.

\section{Center-Unstable Manifold Approximations from GK systems}\label{Sec_GK_centerman}
 We survey in this section classical techniques of approximations of center-unstable manifolds of a steady state for systems of autonomous  ordinary differential equations (ODEs) in $\mathbb{R}^N$, 
for which the right-hand side (RHS) is the RHS of a GK system such as given by \eqref{Galerkin_cptForm}. 

\subsection{The leading-order approximation theorem}
Invariant manifold theory allows for the rigorous derivation of low-dimensional surrogate systems from which not only the system's qualitative behavior near e.g.~a steady state is preserved, but also quantitative features of the nonlinear dynamics are reasonably well approximated such as the solution's amplitude or possible dominant periods.

As a common practice in invariant manifold theory and in view of applications considered in Sec.~\ref{Sec_Reduc_ENSO}, we work with the perturbed equation of Eq.~\eqref{Eq_DDE} (such as Eq.~\eqref{Eq_SS_perturb} below) about a steady state, and its GK approximation of the form 
\be \label{Eq_GK_tau}
\frac{\d \y}{\d t} = A (\tau)\y + G(\y),
\ee
dropping the dependence on $N$.
We assume that nonlinear mapping, $G\colon \mathbb{R}^N \rightarrow \mathbb{R}^N$, satisfies the following tangency condition 
\be \label{Eq_tangencyG}
G(0) = 0, \qquad \text{ and } \qquad  D G(0) = 0. 
\ee
The reduction formulas presented below, are derived  for the case where $G$ does not depend on $\tau$.
 This choice is made to simplify the notations. The general case where $G$ depends on $\tau$ (see \eqref{Galerkin_cptForm}) leads to the same type of formulas as long as the tangency condition \eqref{Eq_tangencyG} is satisfied for all $\tau$.

Assuming that $G$ is sufficiently smooth, then $G(\y)$ admits the following Taylor expansion for $\y$ near the origin:
\begin{equation} \label{G Taylor}
G(\y)= G_k(\underbrace{\y, \cdots, \y}_{k \text{ times}}) + O(\|\y\|^{k+1}),
\end{equation}
where 
\bea \label{k-linear def}
G_k \colon \underbrace{ \mathbb{R}^N \times \cdots \times \mathbb{R}^N}_{k
\text{ times}} \rightarrow \mathbb{R}^N
\eea
denotes a homogenous polynomial of order $k\ge 2$. Often, $G_k(\y)$ is used hereafter as a compact notation for  $G_k(\y,\,\cdots \, , \y)$.
We label the elements of the spectrum $\sigma(A(\tau))$ of $A(\tau)$ according to the lexicographical order. 
According to this rearrangement, the eigenvalues are labeled by a single positive integer $n$, so that 
\bea  \label{eq:ordering-1}
\sigma(A(\tau)) = \{\lambda_n(\tau), \;   1 \leq n \leq N\}, 
\eea
with, for any $1\le n < n'$, either 
\bea
\Res \lambda_{n}(\tau) > \Res \lambda_{n'}(\tau), 
\eea
or
\bea  \label{eq:ordering-3}
\Res \lambda_{n}(\tau) = \Res  \lambda_{n'}(\tau), \; \text{ and } \; \Ims \lambda_{n}(\tau) \geq \Ims \lambda_{n'}(\tau). 
\eea
In this convention, an eigenvalue of algebraic multiplicity $m_c$, is repeated $m_c$ times.

 We are concerned with describing how linear instabilities translate to the nonlinear dynamics.
To do so, the onset of instability is described in terms of the Principle of Exchange of Stability (PES) \cite{MW05}, concerned 
with the loss of stability of the basic steady state.
More precisely,  the PES describes situations for which the spectrum of $A(\tau)$ experiences the following change at a critical parameter $\tau_c$:
\begin{equation} \label{PES}
\begin{aligned}
& \Res \lambda_j(\tau)
\begin{cases} <0 & \mbox{if } \tau < \tau_c, \\ =0 & \mbox{if } \tau = \tau_c,\\ >0 & \mbox{if } \tau > \tau_c,
\end{cases} &&   1 \leq j \leq  m_c, \\
& \Res \lambda_j(\tau_c) < 0, &&  j\ge m_c+1,
\end{aligned}
\end{equation}
for some $m_c>0$, and for $\tau$ in some neighborhood $\mathcal{U}$ of $\tau_c$. 
 Of course, the PES holds also when one crosses $\tau_c$ from above with $ \Res \lambda_j(\tau)<0$ for $\tau>\tau_c$, while $ \Res \lambda_j(\tau)>0$ when $\tau<\tau_c$.

Whatever the situation (destabilization by crossing from above or below a critical value), one associates to the PES condition the following decomposition of the spectrum $\sigma(A(\tau))$:
\bea \label{splitting}
 & \sigma_{\c}(A(\tau)) = \{\lambda_j(\tau) \:|\: j = 1, \cdots, m_c\},\\  
 &\sigma_{\s}(A(\tau)) = \{\lambda_j(\tau) \:|\: j = m_c+1, \cdots, N\}.
\eea

The PES condition prevents other eigenvalues in $\sigma_{\s}(A(\tau))$ from crossing the imaginary axis as $\tau$ varies in $\mathcal{U}$. Hence, no eigenvalues other than those of $\sigma_\c(A(\tau))$ change sign in  the neighborhood $\mathcal{U}$. Furthermore, the PES condition implies, by a continuity argument,  the following 
uniform spectral gap by possibly reducing $\mathcal{U}$ accordingly \cite[Lemma 6.1]{CLW15_vol1},
\be\label{general gap}
0 > 2 k \eta_{\c} > \eta_{\s}, 
\ee
where $k$ is the leading order of $G$, and  
\beas \label{general etac etas}
& \eta_{\c}= \inf_{\tau \in \mathcal{U}} \inf_{j = 1, \cdots, m_c} \mathrm{Re} (\lambda_j(\tau)), & \eta_{\s}= \sup_{\tau \in \mathcal{U}} \sup_{j \ge m_c+1} \mathrm{Re}(\lambda_j(\tau)).
\eeas

To the modes that lose their stability according to the PES, we associate the reduced state  space, $H_{\c}$,  given by 
\be\label{Eq_Hc}
H_{\c} = \mathrm{span}\{\boldsymbol{e}_1, \cdots, \boldsymbol{e}_{m_c}\},
\ee
while a mode $\boldsymbol{e}_n$ with $n\geq m_c+1$ denotes a stable mode. 
Throughout this article we chose not to make explicit the $\tau$-dependence of the eigenmodes but this dependence should be kept in mind. 
We denote by $H_\s$ the subspace spanned by these stable modes.
The projector onto the subspace $H_\s$ (resp.~$H_\c$) is denoted by $\Pi_\s$ (resp.~$\Pi_\c$). Similarly, $A_\c(\tau)$ (resp.~$A_\s(\tau)$) denotes  $\Pi_\c A(\tau)$ (resp.~$\Pi_\s A(\tau)$), and $\y_c$ (resp.~$\y_\s$) denotes the vector in $H_\c$ (resp.~$H_\s$) of the low-mode amplitudes (resp.~stable-mode amplitudes). The inner product in $\mathbb{C}^N$ is denoted by $\langle \cdot, \cdot \rangle$ and is defined by 
\be\label{Eq_inner_pdct}
\langle \bm{a}, \bm{b} \rangle = \sum_{i= 1}^N a_i b_i^\ast, \;\bm{a}, \bm{b} \in \mathbb{C}^N.
\ee 
In what follows we also denote by $ \boldsymbol{e}_j^\ast $ the eigenmodes of the conjugate transpose $A(\tau)^\ast$ of $A(\tau)$.

 Condition \eqref{general gap} ensures in particular that the following spectral gap holds for $\tau$ in $\mathcal{U}$
\bes
\gamma_{m_c}(\tau) = \mathrm{Re}(\lambda_{m_c}(\tau))-\mathrm{Re}(\lambda_{m_c+1}(\tau))>0.
\ees
It is well known that the existence of a (local) exact parameterization or say in other words, of a local $m_c$-dimensional invariant manifold, is subject to the following {\it spectral gap condition}:
\be\label{Eq_spectral_gapcond}
\gamma_{m_c} (\tau) \geq C \mbox{Lip}(G\vert_{\mathcal{V}}),
\ee
 where $\mbox{Lip}(G\vert_{\mathcal{V}})$ denotes the Lipschitz constant of the nonlinearity $G$ restricted to some neighborhood  $\mathcal{V}$ of the origin in $\mathbb{C}^N$, and $C >0$ is typically independent on $\mathcal{V}$. 
 Due to the tangency condition \eqref{Eq_tangencyG},  condition \eqref{Eq_spectral_gapcond} always holds once $\mathcal{V}$  is chosen sufficiently small.  The theory of local invariant manifolds makes thus sense if solutions with sufficiently small amplitudes lie within the appropriate neighborhood $\mathcal{V}$.
 
 In the context of e.g.~nonlinear oscillations that bifurcate from a steady state, the local invariant manifold provides an exact parameterization\footnote{As provided for instance by a center manifold or the unstable manifold of the origin.} of the stable limit cycle near criticality in the case of a supercritical Hopf bifurcation. In the case of subcritical Hopf bifurcation, it provides the parameterization of the unstable limit cycle that emerges in a continuous fashion from the steady state. In Sec.~\ref{Sec_Reduc_ENSO}, we show that the approximation formulas provided by Theorem \ref{Thm_CM_approx} and in Sec.~\ref{Sec_Higher-order} below may allow for approximating not only such  unstable limit cycles that unfold continuously from the origin but also the stable limit cycles that are distant from the origin, corresponding to a ``jump'' transition \cite{MW14} associated with a subcritical Hopf bifurcation.

\bt \label{Thm_CM_approx}
Assume that $G$ and $A(\tau)$ satisfy the assumptions recalled above, and that the PES condition \eqref{PES} is satisfied.   Then for each $\tau$ in a neighborhood $\mathcal{U}$ of $\tau_c$,  Eq.~\eqref{Eq_GK_tau} admits a local invariant manifold, $\mathfrak{M}_\tau=\mbox{graph}(h_\tau)$, with $h_\tau$ that maps $H_\c$ into the stable subspace $H_\s$ such that $h_\tau(0)=0$ and $Dh_\tau(0)=0$ .

Assume that the following non-resonance condition is satisfied:
\begin{equation}  \label{Eq_NR} \tag{NR}
\begin{aligned} 
& \Forall \, (j_1, \cdots, j_k ) \in (1,\cdots,m_c)^k,  \; n\geq m_c+1,  \text{ it holds that} \\
 & \Bigl (\langle G_k(\boldsymbol{e}_{j_1}, \cdots, \boldsymbol{e}_{j_k}), \boldsymbol{e}_n^\ast \rangle \neq 0 \Bigr) \Longrightarrow  \biggl ( \Res (\sum_{\ell=1}^{k} \lambda_{j_\ell} - \lambda_n) \neq 0 \biggr),  
\end{aligned}
\end{equation}
where $G_k$ denotes the leading-order term in the Taylor expansion of $G$. 

Then, the Lyapunov-Perron integral, $\Phi_\tau(\X)= \int_{-\infty}^0 e^{- s A_\s(\tau) } \Pi_{\s} G_k(
e^{s A_\c(\tau)} \X) \, \mathrm{d}s$, is well defined and is a solution to the following homological equation
 \be\label{Eq_homoligical}
\mathcal{L}_{A(\tau)}[\psi] (\X) = \Pi_{\s}G_k(\X),
\ee
where  $\mathcal{L}_A$ denotes the differential operator acting on smooth mappings $\psi$ from $H_\c$ into $H_\s$, defined as:
\be \label{Def_L}
\mathcal{L}_{A(\tau)}[\psi] (\X)=D \psi(\X) A_{\c}(\tau) \X  - A_{\s}(\tau) \psi(\X), \; \X\in H_\c.
\ee

Moreover, $\Phi_\tau(\X)$ provides  the leading-order approximation of the invariant
manifold function $h_\tau$ characterizing $\mathfrak{M}_\tau$, in the sense that  
\be \label{Eq_approx_error}
 \norm{h_\tau(\X) - \Phi_\tau(\X)}=o(\norm{\X}^k), 
\ee
as long as $\X$ lies in a neighborhood $\mathcal{N}(\tau)$ of the origin in $H_\c$ spanned  by the $m_c$ eigenmodes $\boldsymbol{e}_j$ losing their stability (see \eqref{Eq_Hc}), as $\tau$ crosses $\tau_c$.

Finally, $\Phi_\tau(X)$ possesses the following explicit formula:
\be
\label{Eq_phi_leadingorder}
\Phi_\tau(X)=\sum_{n\geq m_c+1}  \Phi_{n,\tau}(\X) \boldsymbol{e}_n, \; \qquad \X\in H_\c,
\ee
where 
\be \label{Phi_n_general}
 \Phi_{n,\tau}(\X) = \sum_{1 \leq j_1, \cdots, j_k \leq m_c}  G_{j_1\cdots j_k}^n M_{j_1\cdots j_k}^{n, \tau} X_{j_1} \cdots X_{j_k}, \; \qquad X_j=  \langle \X, \bm{e}_j \rangle, 
\ee
with 
 \be \label{M_eq} 
M_{j_1\cdots j_k}^{n, \tau}=\Big(\sum_{\ell=1}^{k}\lambda_{j_\ell}(\tau)-\lambda_n(\tau)\Big)^{-1},
\ee
 and the $G_{j_1\cdots j_k}^n$ denoting the  coefficients accounting for the magnitude carried out by $ \boldsymbol{e}_n^\ast $, of the nonlinear interactions through the leading-order term $G_k$, between the low modes $\boldsymbol{e}_{j_1}$, $\cdots$, $\boldsymbol{e}_{j_k}$ (in $H_c$), namely:
\be
G_{j_1\cdots j_k}^n=\langle G_k(\boldsymbol{e}_{j_1}, \cdots, \boldsymbol{e}_{j_k}), \boldsymbol{e}_n^\ast \rangle, \;\; 1\leq j_1, \cdots, j_k \leq m_c.
\ee
\et

Conditions similar to \eqref{Eq_NR} arise in the smooth linearization of dynamical systems near an equilibrium  \cite{sell1985smooth}. Here, condition \eqref{Eq_NR} implies that the eigenvalues of the stable part satisfy a Sternberg condition of order $k$ \cite{sell1985smooth} with respect to the eigenvalues associated with the modes spanning the reduced state space $H_\c$.

 This theorem is essentially a consequence of  \cite[Theorems 1 and 2]{CLM19_closure}, in which condition \eqref{Eq_NR} is a stronger version of that  used for \cite[Theorem 2]{CLM19_closure}; see also \cite[Remark 1 (iv)]{CLM19_closure}. This condition is necessary and sufficient here for $\int_{-\infty}^0 e^{- s A_\s(\tau) } \Pi_{\s} G_k(
e^{s A_\c(\tau)} \X) \, \mathrm{d}s$ to be well defined.  For the derivation of \eqref{Eq_homoligical} see that of Eq.~(4.6) in  \cite{CLM19_closure}. The error estimate \eqref{Eq_approx_error} is a corollary of the more general results  \cite[Theorem 6.1]{CLW15_vol1} and \cite[Corollary 6.1]{CLW15_vol1} proved for stochastic Partial Differential Equations (PDEs). See also \cite[Thm.~6.2.3]{Hen81} and \cite[Lemma 6.2.4]{Hen81}.

The assumptions of Theorem \ref{Thm_CM_approx} ensure that for any $\tau$ in $\mathcal{U}$, the following reduction principle  holds:
\bi
\item[(i)] Any solution $\y(t)$ to Eq.~\eqref{Eq_GK_tau} such that $\y(t_0)$ belongs to $\mathfrak{M}_\tau$ for some $t_0$, stays on   $\mathfrak{M}_\tau$ over an interval of time $[t_0,t_0+\sigma)$, $\sigma>0$, i.e.
\be\label{Eq_inv_local}
\y(t)=\y_\c(t)+h_\tau(\y_\c(t)), \; t\in [t_0,t_0+\sigma),
\ee
where $\y_\c(t)$ denotes the projection of $\y(t)$ onto the subspace $H_{\c}$. 
\ei
\bi
\item[(ii)] If there exists a trajectory $t\mapsto \y(t)$ such that $\y_\c(t)$ belongs, for all $-\infty<t<\infty$, to the neighborhood $ \mathfrak{B}(\tau)$ (in $H_\c$) over which $\mathfrak{M}_\tau$ is well defined,   then the trajectory must lie on $\mathfrak{M}_\tau$.
\ei
Such a formulation of the reduction principle can be inferred for instance from \cite[Theorem 71.4]{SY02} and \cite{Crawford91}.  Property (ii) implies that an invariant set $\Sigma$ to Eq.~\eqref{Eq_GK_tau} of any type, e.g., equilibria, periodic orbits, invariant tori,
must lie in $\mathfrak{M}_\tau$ if its projection onto $H_\c$ is contained in  $\mathfrak{B}(\tau)$, i.e.~if $\Pi_\c \Sigma \subset \mathfrak{B}(\tau)$. Property \eqref{Eq_inv_local} holds then globally in time for the solutions that composed such invariant sets, and thus the knowledge of the $m_c$-dimensional variable, $\y_\c(t)$, is sufficient to entirely determine any solution $\y(t)$ that belongs to such an invariant set; see \cite[Theorem 71.4]{SY02}. Furthermore, $\y_\c(t)$ is obtained as the solution of the following reduced $m_c$-dimensional problem
\be\label{Eq_reduced_absract}
\dot{\widehat{\y}}= A_\c (\tau) \widehat{\y} + \Pi_\c G(\widehat{\y}+ h_\tau(\widehat{\y})), \qquad \widehat{\y}(0)=\y_\c(0) \in \mathfrak{B}(\tau),
\ee
which in turn characterizes the solution $\y(t)$ in $\Sigma$, since the slaving relationship $\y_\s(t)=h_\tau(\y_\c(t))$ holds for any solution $\y(t)$ that belongs to an invariant set $\Sigma$ for which $\Pi_\c \Sigma \subset \mathfrak{B}(\tau)$.

Based on the approximation Theorem \ref{Thm_CM_approx} and the reduction principle, it is thus reasonable to anticipate that the following effective reduced equation,
\be\label{Eq_reduced_applied}
\dot{\x}= A_\c (\tau)  \x + \Pi_\c G_k(\x + \Phi_\tau(\x)), \qquad \x(0)=\y_\c(0) \in \mathfrak{B}(\tau),
\ee
provides an approximation of the invariant set $\Sigma$ for which $\Pi_\c \Sigma \subset\mathfrak{B}(\tau)$.
Since Eq.~\eqref{Eq_reduced_applied} is derived in practice from a GK system, it will be referred to as an {\it effective reduced GK system} or {\it mD reduced GK system}, where $m=\mbox{dim}(H_\c)$. Eq.~\eqref{Eq_reduced_applied} is thus particularly suited for effective approximations of invariant sets  to Eq.~\eqref{Eq_GK_tau}, made of solutions of sufficiently small amplitudes, and  that bifurcate  when one crosses a critical parameter $\tau_c$ that destabilizes an equilibrium point. 
Due to the rigorous convergence results as $N\rightarrow \infty$  to DDE solutions from solutions to GK systems \cite{CGLW16}, it is expected that the invariant set made of bifurcating solutions to the original DDE \eqref{Eq_DDE}, can also be well approximated from the solutions of the effective reduced system \eqref{Eq_reduced_applied}, when $N$ is sufficiently large.      
In some applications, one may need though to go to higher-order approximations than the leading-order one in order to obtain more accurate approximations of $\widehat{\y}(t)$ (and thus of $\y(t)$). This point is described next.  Numerical evidences of such efficient approximations are then presented in Sec.~\ref{Sec_Reduc_ENSO} below; see e.g.~Fig.~\ref{Fig_Hopf_SS}.

\br
 In certain applications, we can encounter that, although designed for $\tau$ close to a critical parameter value $\tau_c$, the effective reduced system \eqref{Eq_reduced_applied} may show skills in approximating bifurcating solutions  to the DDE with not necessarily small amplitudes (order 1), as  $\tau$ is crossing other critical values. Sections \ref{Sec_Bif_SNO} and \ref{Sec_approx_orbits} below provide numerical illustrations of such situations in the case of a delay model of  ENSO. 
\er

\subsection{Analytic formulas for higher-order approximations} \label{Sec_Higher-order}
We provide here a simple derivation of higher-order approximations of an invariant manifold. 
The {\it invariance equation} is insightful in that respect.  It is the equation satisfied by the local invariant manifold, $h_\tau$,  to Eq.~\eqref{Eq_GK_tau} in a neighborhood of the origin,  namely
\be \label{eq:invariance}
D  h(\X) [A_{\c} \X + \Pi_{\c} G(\X + h(\X))] -  A_{\s} h(\X) = \Pi_{\s} G(\X + h(\X)),
\ee
in which we have dropped the dependence on $\tau$ to simplify the notations; see \cite[Corollary 6.2.2]{Hen81}.

This functional equation is a nonlinear system of  first order PDEs that cannot be solved in closed form except in special cases. Many methods exist to seek for approximate solutions to Eq.~\eqref{eq:invariance}; see e.g.~\cite{BK98,EvP04,haro2016parameterization}.
Here, we adopt a standard use of power series expansions that we blend with energy content arguments to neglect certain terms and thus simplify the expressions. Indeed, instead of keeping all the monomials at a given degree arising from such an expansion, we filter out terms that carries significantly less energy compared with those that are kept. This elimination procedure relies on the observation that,  close to criticality, the projected ODE dynamics onto the resolved subspace $H_{\c}$ contains typically a large fraction of the solution's total energy. We illustrate this idea to the case where $G(\y) = G_2(\y,\y) + G_3(\y,\y,\y)$ and a cubic approximation is sought.
The higher-order cases proceed in a same fashion and are omitted here for the sake of concision.   
We refer to \cite{haro2016parameterization}  for mathematical details about higher-order approximations of invariant manifolds.

To determine the third-order approximation, we replace $h$ in the invariance equation \eqref{eq:invariance} by $\Phi_{\tau} + \psi$, where $\psi$ denotes the homogeneous cubic terms in the power expansion of $h$, to be determined. By identifying the terms of order two, we recover Eq.~\eqref{Eq_homoligical} satisfied by $\Phi_{\tau}$.  By identifying the terms of order three, we obtain the following equation for $\psi$: 
\be \label{Eq_invariance_psi}
\mathcal{L}_{A}[\psi] (\X) =  -D \Phi_{\tau}(\X) \Pi_{\c} G_2(\X) + \Pi_{\s} G_2(\X, \Phi_{\tau}(\X)) + \Pi_{\s} G_2(\Phi_{\tau}(\X), \X) + \Pi_{\s} G_3(\X), 
\ee
where $\mathcal{L}_{A}[\psi] (\X) = D  \psi(\X) A_{\c} \X  - A_{\s} \psi(\X)$; see \eqref{Def_L}.

Note that the RHS  is a homogeneous cubic polynomial in the $\X$-variable. As recalled above, close to criticality, a large fraction of the energy  is contained in the low modes and therefore the energy carried by $\y_\s$ is typically much smaller than $\|\y_\c\|^2$. It is then reasonable to expect that the energy carried by $\Phi_{\tau}(\X)$ is much smaller than $\|\X\|^2$ for $\X = \y_\c$. Thus, one may assume that in the RHS of \eqref{Eq_invariance_psi}, the term $\Pi_{\s} G_3(\X)$ dominates the other three terms provided that $\|G_2(\y)\|/\|\y\|^2$ is on the same order of magnitude as $\|G_3(\y)\|/\|\y\|^3$. As a consequence, it is reasonable to hope  for reasonably accurate  approximations of  $\psi$  with $h_3$ by  simply solving the equation: 
\be
D  h_3(\X) A_{\c} \X  - A_{\s} h_3(\X) =  \Pi_{\s} G_3(\X). 
\ee
Note that this equation is exactly Eq.~\eqref{Eq_homoligical} with $k=3$. In virtue of Theorem \ref{Thm_CM_approx}, the existence of $h_3$ is guaranteed by non-resonance condition \eqref{Eq_NR} (with $k=3$), and $h_3$ is thus given by \eqref{Eq_phi_leadingorder}--\eqref{Phi_n_general} with $k=3$. We arrive then at the following ``high-mode'' parameterization
\bea \label{Eq_Phi}
\Psi_\tau(\X)&= \Phi_\tau(\X) + h_3(\X)= \sum_{n=m_c+1}^N \Psi_{n,\tau}(X)\bm{e}_n \\
 & = \sum_{n=m_c+1}^N \bigg( \sum_{1\leq j_1, j_2 \leq m_c} \frac{\langle G_2(\boldsymbol{e}_{j_1}, \boldsymbol{e}_{j_2}), \boldsymbol{e}_n^\ast \rangle}{\lambda_{j_1}(\tau) + \lambda_{j_2}(\tau) - \lambda_n(\tau)} X_{j_1} X_{j_2} 
\\ & \hspace{2cm}+ \hspace{-.5cm}\sum_{1\leq j_1, j_2, j_3 \leq m_c} \frac{\langle G_3(\boldsymbol{e}_{j_1}, \boldsymbol{e}_{j_2}, \boldsymbol{e}_{j_3}), \boldsymbol{e}_n^\ast \rangle}{\lambda_{j_1}(\tau) + \lambda_{j_2}(\tau) + \lambda_{j_3}(\tau) - \lambda_n(\tau)} X_{j_1} X_{j_2}X_{j_3}
\bigg)\boldsymbol{e}_n.
\eea
In Sec.~\ref{Sec_Reduc_ENSO} below, we show that such parameterizations enable us to reach accurate approximations of bifurcating solutions in the case of a DDE with cubic nonlinearity. 

\section{Bifurcations Analysis of a Delay El Ni\~no-Southern Oscillation (ENSO) Model}\label{Sec_Reduc_ENSO}

\subsection{GK  approximations of the Suarez and Schopf ENSO model}
The Suarez and Schopf model  is given by the following DDE  \cite{Suarez_al88} 
\be \label{eq:SS}
\frac{\d T}{\d t} = T(t) - \alpha T(t - \tau) - T^3(t),
\ee
where the unknown $T$ represents the non-dimensionalized sea surface temperature (SST) anomalies at the eastern equatorial Pacific Ocean, $\tau$ and $\alpha$ are positive constants, and the physically relevant range of $\alpha$ used in \cite{Suarez_al88} is $(0,1)$. 
For this given range of $\alpha$, Eq.~\eqref{eq:SS} admits three fixed points: 
\be\label{Eq_steady_states}
T_0 = 0, \qquad T_{+} = \sqrt{1- \alpha}, \qquad T_{-} = - \sqrt{1- \alpha}.
\ee
 The first term in Eq.~\eqref{eq:SS} reflects an instantaneous positive feedback, whereby an SST perturbation heats the atmosphere, whose wind response drives ocean currents to reinforce the
original perturbation. The second term in Eq.~\eqref{eq:SS} accounts for the transit time of equatorially trapped oceanic waves to cross the Pacific ocean \cite{boutle2007nino}. The third term accounts then for  saturation effects which limit the instability growth due to the
positive feedback by effects tied to advective processes in the ocean and moist processes in the atmosphere.  A linear stability analysis of this simple model was conducted
in \cite{Suarez_al88}  to show that the steady-state solution loses stability for certain parameter values. The resulting
periodic solutions have periods of at least twice the length of
delay,  supporting that this simple feedback mechanism can be, in theory, consistent with ENSO’s oscillatory behaviour on an interannual timescale.  In Sec.~\ref{Sec_ENSO_var} below, we analyze in greater details the nonlinear structures at play that organize the time-variability of solutions to the Suarez and Schopf  model when subject to noise disturbances. These nonlinear structures are invariant sets that emerge through local and global bifurcations.    

In what follows, we conduct thus a bifurcation analysis of the Suarez and Schopf  model. To do so, we apply the  GK approximation framework of Sec.~\ref{GK_section} and the center-unstable manifold reduction formulas of Sec.~\ref{Sec_GK_centerman}.  In that respect, to fit within the framework of Sec.~\ref{Sec_GK_centerman}, we first derive from Eq.~\eqref{eq:SS},  the DDE satisfied by the perturbed variable about the steady state $T_{+}$, 
\be\label{Eq_u}
\T= T - T_{+},
\ee 
namely,
\be \label{Eq_SS_perturb}
\frac{\d \T}{\d t} = (1 - 3 T_{+}^2) \T(t)- \alpha \T(t - \tau)  - 3 T_{+} \T^2(t) - \T^3(t).
\ee
The above equation fits into Eq.~\eqref{Eq_DDE} with 
\bea \label{eq:SS_pars}
a = 1 - 3 T_{+}^2, \quad b = -\alpha,  \quad c = 0, \text{and} \quad F(\T(t)) = - 3 T_{+} \T^2(t) - \T^3(t).
\eea 
By applying Eqns.~\eqref{Galerkin_cptForm}--\eqref{eq:G} to Eq.~\eqref{Eq_SS_perturb}, we obtain the following $N$-dimensional GK system
\be \label{Eq_Galerkin_SS}
\frac{\d \y}{\d t} = A(\tau) \y  +  G (\y), \qquad \y= (y_1, \cdots, y_N)^{T},
\ee
where the entries of the matrix $A(\tau)$ are given here by
\bea \label{eq:A_SS}
(A(\tau))_{j,n} &  = \frac{1}{\|\mathcal{K}_j\|_{\mathcal{E}}^2 } \sum_{n=0}^{N-1}  \Bigl( 1 - 3 T_{+}^2  - \alpha K_n(-1)  \\
&  \qquad +  \frac{2}{\tau}\sum_{k=0}^{n-1} a_{n,k} \left( \delta_{j,k} \|\mathcal{K}_j\|^2_{\mathcal{E}} - 1 \right ) \Bigr),  \qquad \mbox{for $j, n  = 0, \cdots, N-1$.}
\eea
As in Sec.~\ref{Sec_GKaction}, the $ a_{n,k}$ are determined thanks to Proposition \ref{prop:dPn} recalled in Appendix~\ref{sect:coef_matrix_proof}, and  the $\|\mathcal{K}_j\|_{\mathcal{E}}$ and $K_n(-1)$ are given by \eqref{eq:Kn_norm} and \eqref{Eq_K-1}, respectively.

The nonlinearity $G$ is given as the following sum of $N$-dimensional mappings
\bes 
 G(\y) =  G_2(\y) +  G_3(\y),
\ees
with 
\be \label{Eq_G_decomp_SS}
G_2(\y) = - 3 T_{+}  \left( \sum_{n=0}^{N-1} y_n \right)^2 \boldsymbol{\nu}_N, \mbox{ and } \; G_3(\y) = - \left( \sum_{n=0}^{N-1} y_n \right)^3 \boldsymbol{\nu}_N,
\ee
 where $\boldsymbol{\nu}_N$ denotes the $N$-dimensional column vector  
 \be\label{Eq_vector_norms}
 \boldsymbol{\nu}_N=\Big(\frac{1}{\|\mathcal{K}_0\|_{\mathcal{E}}^2}, \cdots, \frac{1}{\|\mathcal{K}_{N-1}\|_{\mathcal{E}}^2} \Big)^T
\ee
We arrive finally at the following $N$-dimensional GK approximation for the DDE \eqref{Eq_SS_perturb}
\bea  \label{Eq_Galerkin_SS_6D}
\frac{\d \y}{\d t} = \left(\frac{2}{\tau} \mathcal{P}_N +\mathcal{Q}_N\right) \y  - \left(\sum_{n=1}^{N} y_n \right)^2  \left( 3 T_{+}  + \sum_{n=1}^{N} y_n \right) \boldsymbol{\nu}_N.
\eea
We refer to Appendix~\ref{sect:coef_matrix_proof} for the computation of  $\mathcal{P}_N$ and $\mathcal{Q}_N.$
 By construction, the nonlinear term $G$ satisfies the tangency condition \eqref{Eq_tangencyG} and thus one can unpack the center-unstable manifold framework of Sec.~\ref{Sec_GK_centerman}. This is done in Sec.~\ref{Sec_Bif_Hopf} below. To prepare for the resulting bifurcation analysis, we consider the case $\alpha = 0.75$. For this value of $\alpha$, the steady state is given by $T_{+} = \sqrt{1 - \alpha} = 0.5$.

 We conclude this section with Table \ref{Table_Error} that shows the $L^\infty$-error achieved by various $N$-dimensional GK systems in approximating stable limit cycles to the DDE \eqref{Eq_SS_perturb} for $\tau$-values that will be considered later on.

\begin{table}[h] 
\caption{{\bf $L^\infty$-error achieved by the $N$-dimensional GK approximations to  the DDE \eqref{Eq_SS_perturb}.} Here, the $L^\infty$-errors are computed over one  period of the stable limit cycle for each tested $\tau$-value, with the DDE \eqref{Eq_SS_perturb} initialized using a segment on the stable limit cycle, and the GK systems initialized using the projection of the DDE initial data onto the first $N$ Koornwinder basis functions. The timestep size is set to $\delta t = \tau/2^{18}$.}
\label{Table_Error}
\centering
\begin{tabular}{ccccc}
\toprule\noalign{\smallskip}
   GK dimension  & $\tau = 1.562$   &  $\tau = 1.6$ &  $\tau = 1.7$  & $\tau = 1.9$   \\ 
\noalign{\smallskip}\hline\noalign{\smallskip}
 $N=4$ &     $5.66\times 10^{-2}$ & $3.87\times 10^{-2}$ &  $4.29\times 10^{-2}$  &   $5.71\times 10^{-2}$ \\ 
 $N=6$ &     $1.44\times 10^{-4}$ & $1.80\times 10^{-4}$ &  $7.75\times 10^{-4}$  &   $3.60\times 10^{-3}$ \\ 
  $N=8$ &    $1.54\times 10^{-4}$ & $6.95\times 10^{-5}$ &  $2.13\times 10^{-5}$  &   $4.37\times 10^{-4}$ \\ 
   $N=10$ & $1.31\times 10^{-4}$ & $5.60\times 10^{-5}$ &  $3.65\times 10^{-5}$  &   $8.13\times 10^{-5}$ \\  
\noalign{\smallskip} \bottomrule 
\end{tabular}
\end{table}

\subsection{Subcritical Hopf Bifurcation}\label{Sec_Bif_Hopf}
In this section, we apply the center-unstable manifold reduction framework of Sec.~\ref{Sec_GK_centerman} to the the GK system \eqref{Eq_Galerkin_SS_6D}, and in particular  the high-mode parameterization $\Psi_\tau$ given by \eqref{Eq_Phi}, to derive our  effective reduced GK system and conduct a bifurcation analysis of the (perturbed) Suarez and Schopf model \eqref{Eq_SS_perturb}.

Thus, our $m_c$-dimensional effective reduced GK system based on $\Psi_\tau$ writes in compact form
\be \label{Eq_EffectiveReduced_ENSO_cpt} 
\dot{\x} = A_{\c} (\tau) \x + \Pi_{\c} G \bigl ( \x +  \Psi_\tau(\x) \bigr),
\ee
for which $m_c$ is determined by analyzing the modes that destabilize once a critical delay parameter $\tau$ is crossed. 

In the eigenbasis coordinate system of  $A (\tau)$, the GK system \eqref{Eq_Galerkin_SS_6D} writes 
\be \label{Eq_GK_ENSO}
\dot{y}_j = \lambda_j (\tau) y_j + \Bigl \langle G \Bigl( \sum_{\ell = 1}^N y_ \ell \boldsymbol{e}_\ell  \Bigr), \boldsymbol{e}_j^\ast \Bigr \rangle, \;\; j = 1, \cdots, N,
\ee
where $y_j = \langle \y,  \boldsymbol{e}_j^*\rangle$, with inner product defined in \eqref{Eq_inner_pdct}.

\begin{figure}[hbtp]
   \centering
\includegraphics[width=0.6\textwidth, height=0.4\textwidth]{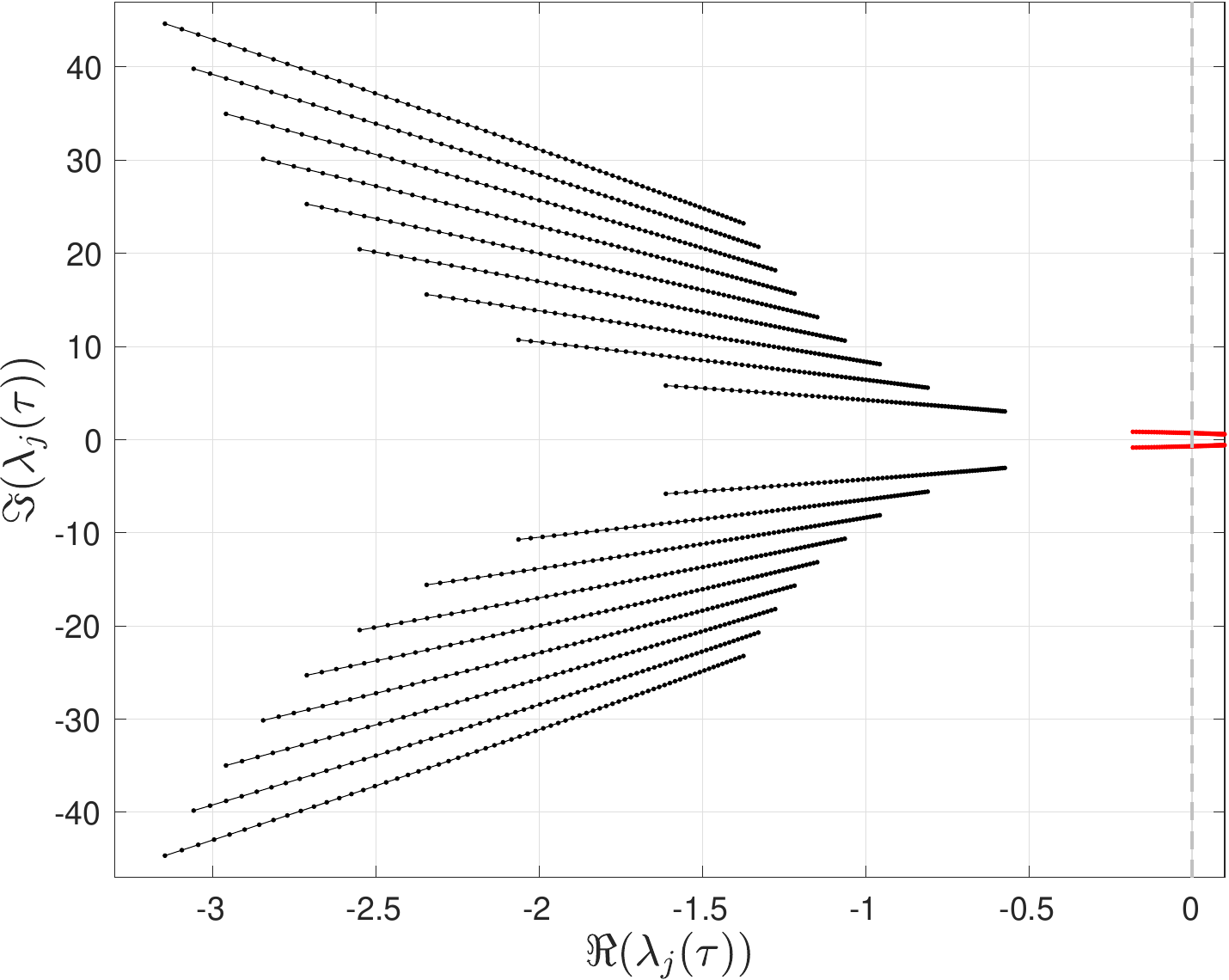}
  \caption{{\bf Eigenvalues dependence as $\tau$ crosses its critical value $\tau_c$}. 
  Are shown here, the first 10 pairs $(\lambda_j(\tau),\overline{\lambda_j(\tau)})$ as $\tau$ is increased  from $\tau=1.3$ to $\tau=2.5$, for $\alpha=0.75$. 
  These pairs move all from left to right.  
 The pair of eigenvalues that crosses the imaginary axis for the critical delay parameter $\tau=\tau_c$ is shown in red, while the stable pairs are shown in black. These eigenvalues are computed from the GK linear part $A(\tau)$ given in \eqref{eq:A_SS} in dimension $N=50$. It is noteworthy that these 10 pairs of GK eigenvalues (for $N=50$) satisfy the  actual characteristic equation associated with the linear part of the DDE \eqref{Eq_SS_perturb}, $\lambda = (1 - 3 T_{+}^2) - \alpha e^{-\lambda \tau}$, up to a maximal  error of $10^{-4}$, as $\tau$ varies from $\tau=1.3$ to $\tau=2.5$.}\label{Fig_eigen}
\end{figure}

We observe that a dominant pair of 
complex eigenvalues crosses the imaginary axis as $\tau$ crosses from below the critical value $\tau_\c \approx  1.74$ (see Table \ref{Table_Lyap_coeff}), whereas the other pairs of eigenvalues stay within the left half complex plane, with a clear spectral gap, even for $\tau$ above the critical value; see Fig.~\ref{Fig_eigen}. We choose thus $m_c=2$ and the subspaces $H_{\c}$ and $H_{\s}$ as follows: 
\be
H_{\c} = \mathrm{span}\{\boldsymbol{e}_1, \boldsymbol{e}_2\}, \qquad H_{\s} = \mathrm{span}\{\boldsymbol{e}_{3}, \cdots, \boldsymbol{e}_N\},
\ee
in which $(\bm{e}_1,\bm{e}_2)$ denotes the conjugate pair of modes that destabilize as $\tau$ crosses $\tau_c$. 
Setting $m_c=2$ in \eqref{Eq_Phi}, the 2D reduced GK system \eqref{Eq_EffectiveReduced_ENSO_cpt} (in $\mathbb{C}^2$)  writes in the eigenbasis coordinates:  
\bea \label{Eq_EffectiveReduced_ENSO}
\dot{x}_j= \lambda_j(\tau) x_j + \langle G( \x + \Psi_{\tau}(\x)), \boldsymbol{e}_j^\ast \rangle, \qquad j= 1, 2,
\eea
with $\x=x_1 \boldsymbol{e}_1 + x_2 \boldsymbol{e}_2$, and where the parameterization $\Psi_{\tau}$ of the neglected modes in $H_{\s}$ is given by \eqref{Eq_Phi}.

The behavior of the eigenvalues shown in Fig.~\ref{Fig_eigen} is an illustration of the PES condition \eqref{PES} for the $N$-dimensional GK approximation~\eqref{Eq_GK_ENSO} of the DDE \eqref{Eq_SS_perturb}. In particular, the critical crossing of the dominant pair (red curves in Fig.~\ref{Fig_eigen}) from the left to the right half plane  is a manifestation of the well-known transversality condition, for the underlying GK system to admit a Hopf bifurcation \cite[Eqns.~(32)-(33)]{CKL20}.  Theorem III.1 of \cite{CKL20} provides then the precise conditions for a Hopf bifurcation to take place for the GK system and ensures that its type,  either supercritical or subcritical, is determined by calculating the Lyapunov coefficient from the 2D reduced GK system \eqref{Eq_EffectiveReduced_ENSO} put into its Stuart-Landau (SL) form \cite[Eq.~(36)]{CKL20}.   
As examined in \cite{CKL20}, the same type of Hopf bifurcation such as predicted by the SL equation, occurs then for the DDE as well,  as long as the coefficients of the SL equation are calculated from a sufficiently large dimensional GK system.

 In that respect, we determine the Lyapunov coefficient as given by the analytic formula \cite[Eq.~(40)]{CKL20} from the coefficients of the $N$-dimensional GK approximation~\eqref{Eq_GK_ENSO} of the DDE \eqref{Eq_SS_perturb}.
To this end,  we rewrite $G_2$ and $G_3$ given in \eqref{Eq_G_decomp_SS} as the following bilinear and trilinear forms: 
\bea
G_N^{(2)}(\boldsymbol{p}, \boldsymbol{q}) & = - 3 T_{+}  \left( \sum_{n=0}^{N-1} p_n \right) \left( \sum_{n=0}^{N-1} q_n \right) \boldsymbol{\nu}_N \\
G_N^{(3)}(\boldsymbol{p}, \boldsymbol{q},\boldsymbol{r}) & = - \left( \sum_{n=0}^{N-1} p_n \right)\left( \sum_{n=0}^{N-1} q_n \right)\left( \sum_{n=0}^{N-1} r_n \right) \boldsymbol{\nu}_N.
\eea
The Lyapunov coefficient $\ell^N_1(\tau_c)$ is then given by 
\be  \label{Eq_l1_GK}
\ell_1^N (\tau_c)=  \mathrm{Re}\Bigl (\frac{ a^N_{20}a^N_{11} \sqrt{-1} }{\mathrm{Im}(\lambda_1(\tau_c))}  + a^N_{21} \Bigr),
\ee
where 
\bea \label{eq:l1_coef1}
a_{20}^N & = \langle G_N^{(2)}(\boldsymbol{e}_1, \boldsymbol{e}_1), \boldsymbol{e}_1^\ast \rangle,\\
a^N_{11}  &= \langle G_N^{(2)}(\boldsymbol{e}_1, \boldsymbol{e}_2), \boldsymbol{e}_1^\ast \rangle + \langle G_N^{(2)}(\boldsymbol{e}_2, \boldsymbol{e}_1), \boldsymbol{e}_1^\ast \rangle, 
\eea
and
\bea \label{eq:l1_coef2}
a_{21}^N& = \langle G_N^{(3)}(\boldsymbol{e}_1, \boldsymbol{e}_1, \boldsymbol{e}_2), \boldsymbol{e}_1^\ast \rangle + \langle G_N^{(3)}(\boldsymbol{e}_1, \boldsymbol{e}_2, \boldsymbol{e}_1), \boldsymbol{e}_1^\ast \rangle + \langle G_N^{(3)}(\boldsymbol{e}_2, \boldsymbol{e}_1, \boldsymbol{e}_1), \boldsymbol{e}_1^\ast \rangle \\
&  \qquad + \sum_{n=3}^{Nd} \frac{\langle G_N^{(2)}(\boldsymbol{e}_{1}, \boldsymbol{e}_{2}), \boldsymbol{e}_n^\ast \rangle + \langle G_N^{(2)}(\boldsymbol{e}_{2}, \boldsymbol{e}_{1}), \boldsymbol{e}_n^\ast \rangle }{ 2  \mathrm{Re} (\lambda_{1}(\tau_c))  - \lambda_n(\tau_c)} \Bigl[ \langle G_N^{(2)}(\boldsymbol{e}_1, \boldsymbol{e}_n), \boldsymbol{e}_1^\ast \rangle + \langle G_N^{(2)}(\boldsymbol{e}_n, \boldsymbol{e}_1), \boldsymbol{e}_1^\ast \rangle \Bigr] \\
&  \qquad + \sum_{n=3}^{Nd} \frac{\langle G_N^{(2)}(\boldsymbol{e}_{1}, \boldsymbol{e}_{1}), \boldsymbol{e}_n^\ast \rangle}{ 2 \lambda_{1}(\tau_c) - \lambda_n(\tau_c)} \Bigl[ \langle G_N^{(2)}(\boldsymbol{e}_2, \boldsymbol{e}_n), \boldsymbol{e}_1^\ast \rangle + \langle G_N^{(2)}(\boldsymbol{e}_n,\boldsymbol{e}_2), \boldsymbol{e}_1^\ast \rangle \Bigr].
\eea
Theorem III.1 of \cite{CKL20} ensures then that subcritical Hopf bifurcation occurs for the $N$-dimensional GK system~\eqref{Eq_GK_ENSO} if $\ell^N_1(\tau_c)>0$  and a supercritical Hopf bifurcation occurs  if $\ell^N_1(\tau_c) < 0$. 

 Table~\ref{Table_Lyap_coeff} reports on the $\tau$-value $\tau_c$  at which the dominant pair $(\lambda_1(\tau),\overline{\lambda_1(\tau)})$ crosses the imaginary axis, along with its corresponding Lyapunov coefficient $\ell^N_1(\tau_c)$, in terms of the   GK dimension $N$.
As can be observed, both  $\tau_c$ and $\ell^N_1(\tau_c)$  converges quickly as the GK dimension $N$ increases. 
Since $\ell^N_1(\tau_c)$  converges  to a positive value, we infer that the original Suarez and Schopf model~\eqref{Eq_SS_perturb} admits a subcritical Hopf bifurcation for the parameter regime considered here ($\alpha=0.75$).  

\begin{figure}[hbtp]
   \centering
\includegraphics[width=0.5\textwidth, height=0.3\textwidth]{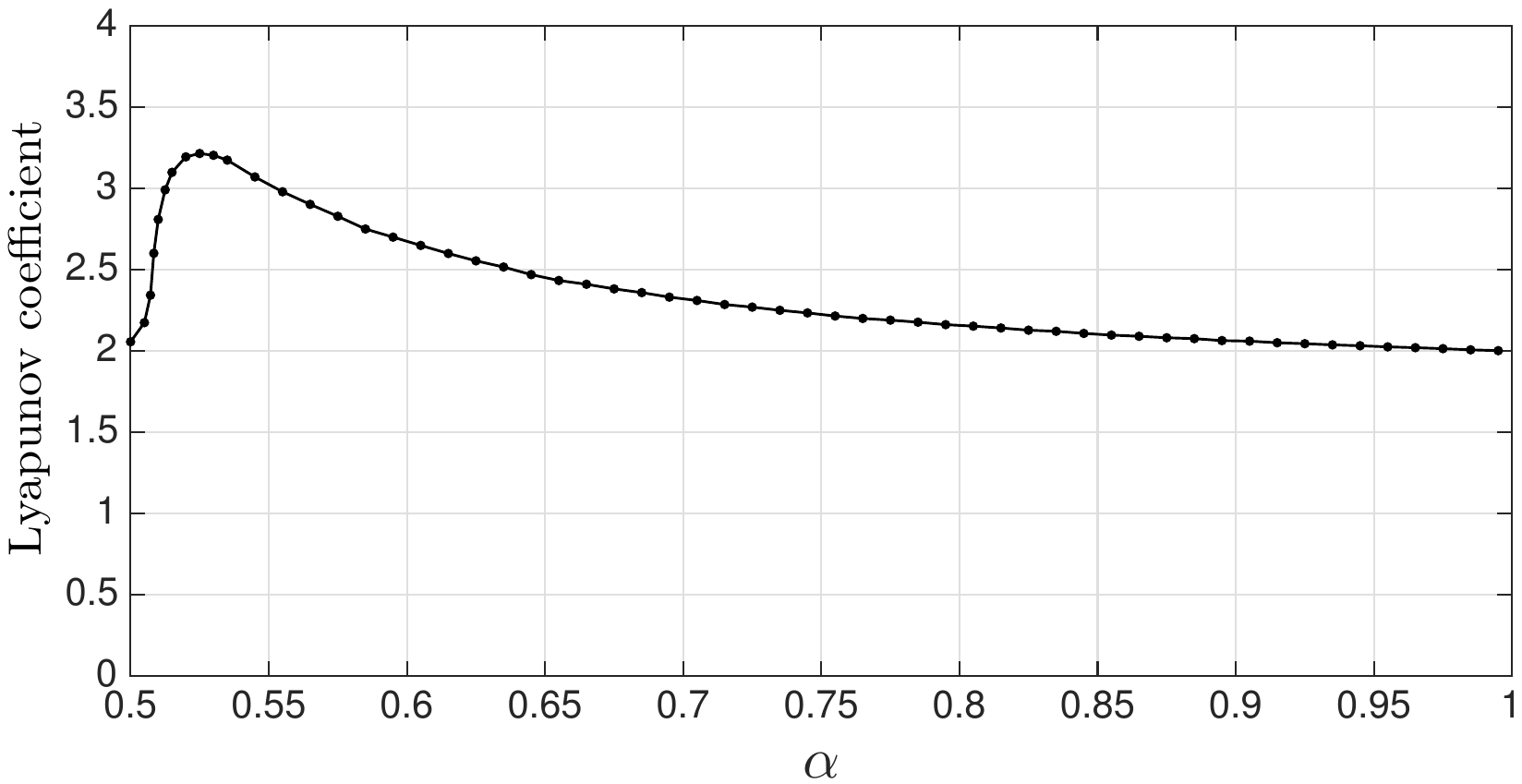}
  \caption{{\bf Lyapunov coefficient $\ell^N_1(\tau_c)$ given in Eq.~\eqref{Eq_l1_GK}: $\alpha$-dependence.}  Here, the critical value $\tau_c$ and the entries in Eq.~\eqref{Eq_l1_GK} are computed from a high-dimensional GK approximation   ($N=20$) of the (perturbed) Suarez and Schopf model \eqref{Eq_SS_perturb}.}\label{Fig_Lyap_coeff}
\end{figure}

 The Laypunov coefficient $\ell^N_1(\tau_c)$ can obviously be computed for other values of $\alpha$ in  $(0.5,1)$, once $\tau_c$ is determined as the critical $\tau$-value at which the dominant pair of eigenvalues  of $A(\tau)$ crosses the imaginary axis as in Fig.~\ref{Fig_eigen}.  The results are shown in Fig.~\ref{Fig_Lyap_coeff}  for $N=20$. The positivity of $\ell^N_1(\tau_c)$ shows that the steady state $T_+$ undergoes always a subcritical Hopf bifurcation for $\alpha$  in $(0.5,1)$. Note that this is the maximum range of $\alpha$ values for which $T_+$ can lose stability, since when $\alpha > 1$, $T_+ = \sqrt{1-\alpha}$ ceases to exist, and when $\alpha < 0.5$, $T_+$ is always linearly stable for any $\tau$.

One might wonder though whether $\tau_c$ as predicted from the GK systems converges to the actual value of  $\tau_c$, obtained by analyzing the  characteristic equation associated with the DDE \eqref{Eq_SS_perturb}. We answer this questioning by the affirmative.  Recall indeed that one can infer the following  analytic expression $\tau^\text{true}_c $ given by \cite[Eq.~(7)]{boutle2007nino} from the characteristic equation,  
\be \label{Eq_tauc_true}
\tau^\text{true}_c =  \frac{\arccos(\frac{3\alpha - 2}{\alpha})}{\sqrt{\alpha^2 - (3\alpha - 2)^2}}.
\ee
For $\alpha = 0.75$ considered here, we have $\tau^\text{true}_c \approx 1.740839502734206$. As shown in Table~\ref{Table_Lyap_coeff}, the $10$-dimensional GK system already provides a highly accurate approximation of $\tau^\text{true}_c$ with precision $10^{-7}$. A similar precision is obtained for other values of $\alpha$ in  $(0.5,1)$.

\begin{table}[h] 
\caption{{\bf The estimated $\tau_c$ and Lyapunov coefficient $\ell^N_1(\tau_c)$ as GK dimension $N$ increases.}}
\label{Table_Lyap_coeff}
\centering
\begin{tabular}{cccccc}
\toprule\noalign{\smallskip}
  &  $N=4$ &  $N=6$  &  $N=8$ & $N=10$ & $N=12$   \\ 
\noalign{\smallskip}\hline\noalign{\smallskip}
$\tau_c$ & 1.7343471   & 1.7408640 &   1.7408394 &   1.7408395 &   1.7408395 \\ 
\noalign{\smallskip}
$\ell^N_1(\tau_c)$ &  2.1884711 &  2.2246906  &   2.2247594  &   2.2247568 &   2.2247568   \\
\noalign{\smallskip} \bottomrule 
\end{tabular}
\end{table}

\subsection{Saddle-Node bifurcation of periodic Orbits (SNO) and homoclinic orbit}\label{Sec_Bif_SNO}
\begin{figure}[hbtp]
   \centering
\includegraphics[width=0.8\textwidth, height=0.5\textwidth]{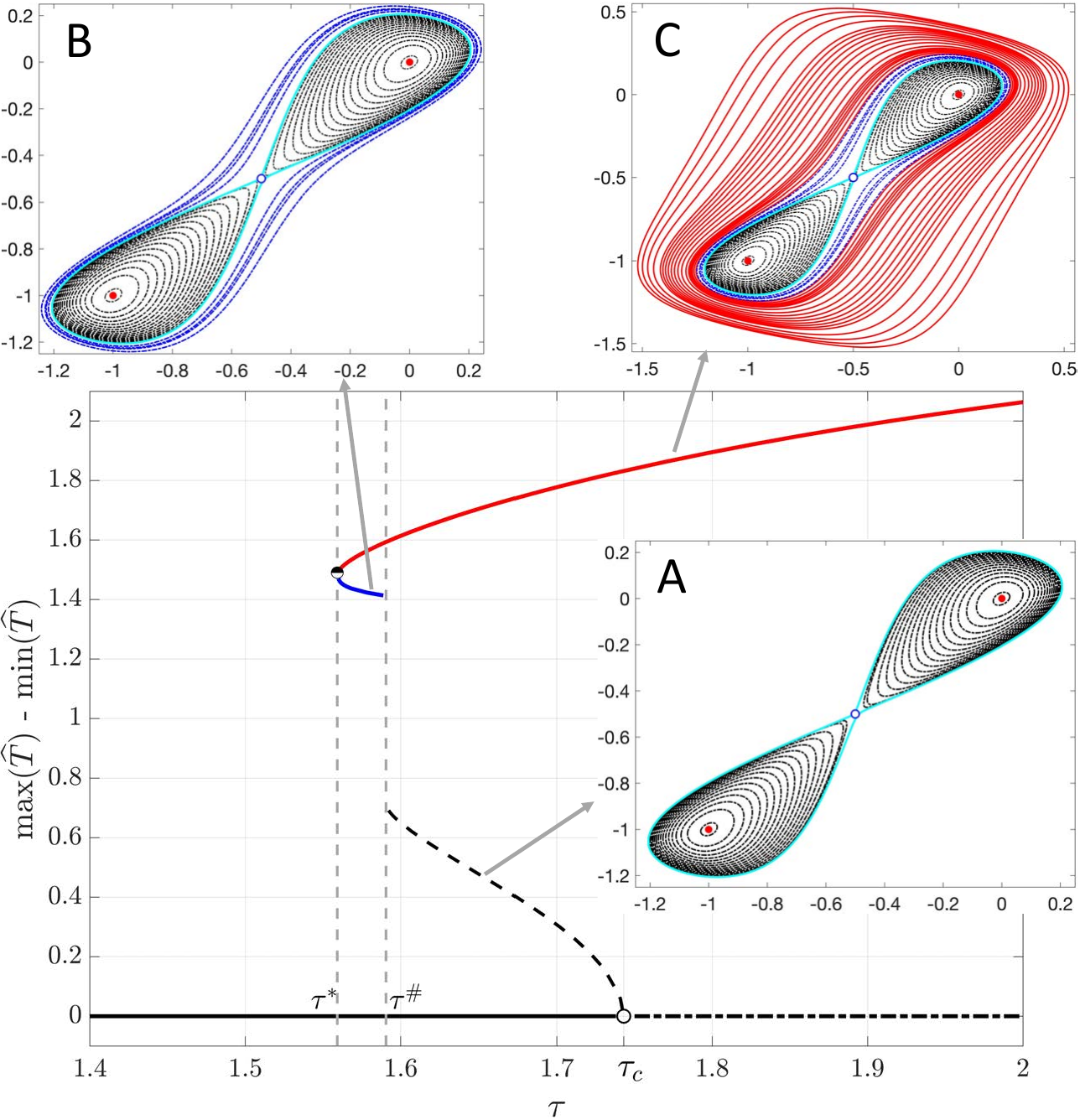}
\caption{{\bf Bifurcation Diagram from the 2D Reduced GK System \eqref{Eq_EffectiveReduced_ENSO}}. This diagram  describes precisely the local and global bifurcations occurring for the (perturbed) Suarez and Schopf model \eqref{Eq_SS_perturb}, given the approximation skills of the reduced system \eqref{Eq_EffectiveReduced_ENSO} in approximating  the DDE's periodic orbits; see Sec.~\ref{Sec_approx_orbits} below. As for Fig.~\ref{Fig_Hopf_SS} below, it turned out to be sufficient to use $N=6$ in the construction of the parameterization $\Psi_{\tau}$ involved in Eq.~\eqref{Eq_EffectiveReduced_ENSO}, to reach such skills. 
In each inset, the stable steady states $T_+$ and $T_-$ given by Eq.~\eqref{Eq_steady_states} are marked by red dots with $T_+$ corresponding to $(0,0)$ since this diagram is computed for the perturbed DDE \eqref{Eq_EffectiveReduced_ENSO}. This diagram reads as follows.
 The steady state $T_+$ (resp.~ $T_-$) is locally stable for all $\tau < \tau_c$, and loses its stability through a subcritical Hopf bifurcation at $\tau_c \approx 1.7408$; see Table~\ref{Table_Lyap_coeff} and Eq.~\eqref{Eq_tauc_true}. As $\tau$ approaches $\tau^{\sharp} \approx 1.5906$ from above, the bifurcating UPOs (black dashed curves) merge into  a homoclinic orbit (cyan curve) connecting the stable and unstable directions of the saddle steady state $T_0$, marked by an empty blue circle in inset A. This merging of UPOs terminates the branch of UPOs shown as dashed black line in the diagram. 
The diminishing of $\tau$ below $\tau^{\sharp}$ leads to new  UPOs that encompass the homoclinic orbit, with amplitude that grows as $\tau$ approaches $\tau^{\sharp}$; see blue dashed curves in inset B. This latter branch of UPOs terminates at $\tau^* \approx 1.5592$ through a SNO bifurcation that gives rise to a branch of stable periodic orbits shown  by the red curves in inset C.}\label{Fig_bif_combo}
\end{figure}

We show in this Section, that the 2D reduced GK system \eqref{Eq_EffectiveReduced_ENSO} is not only useful to predict the subcritical Hopf bifurcation occurring for the Suarez and Schopf model \eqref{Eq_SS_perturb} at the parameter value $\tau=\tau_c$, but also bifurcations taking place at other critical values of $\tau$. For instance, a simple backward integration of  Eq.~\eqref{Eq_EffectiveReduced_ENSO} allows for computing a whole family of unstable periodic orbits (UPOs) unfolding from subcritical Hopf bifurcations from $T_+$ and $T_{-}$ (red points in  Fig.~\ref{Fig_bif_combo}A) as $\tau$ is further decreased to a value $\tau^{\sharp}$ for which the UPOs (dashed black curves in Fig.~\ref{Fig_bif_combo}A) hit a saddle equilibrium from both sides, resulting into a homoclinic orbit when $\tau$ reaches a critical value $\tau^{\sharp}$; see cyan curve in Fig.~\ref{Fig_bif_combo}A.  The presence of this homoclinic orbit is manifested by the jump displayed by the curve of UPO's amplitudes\footnote{Obtained by approximating  $\T(t)$ solving Eq.~\eqref{Eq_SS_perturb} by the UPOs from the reduced equation according to \eqref{Eq_SS_reconstruct2} below.}, $\tau\mapsto \max(\T(t)) - \min(\T(t))$, in the diagram of Fig.~\ref{Fig_bif_combo}.  
This jump is explained by our way of  calculating them as these amplitudes are indeed computed  over $(\tau^{\sharp},\tau_c)$ only from the branching family of UPOs emanating from $T_{+}$, i.e.~located in the lobe encircling $T_+$ of this homoclinic orbit.

After this jump, as $\tau$ is further decreased from  $\tau^{\sharp}$, the UPOs encompass the homoclinic orbit (blue curves in Fig.~\ref{Fig_bif_combo}B), and their amplitude keeps increasing,  until eventually, loosing their instability into a Saddle-Node bifurcation of periodic Orbits (SNO bifurcation) for $\tau=\tau^\ast$. At this critical value the periodic orbit shows mixed stability, with its basin of attraction corresponding to the exterior of the closed curve.  Finally, as $\tau$ is increased from $\tau^\ast$, stable limit cycles unfold with increasing amplitude (red curves in Fig.~\ref{Fig_bif_combo}C) and pronounced nonlinear features expressed by their ovaloid shape.

As explained in the next section, this bifurcation diagram obtained from the  2D reduced GK equation \eqref{Eq_EffectiveReduced_ENSO} describes precisely the local and global bifurcations occurring for 
the (perturbed) Suarez and Schopf model \eqref{Eq_SS_perturb} due to  the remarkable approximation skills shown by the  reduced equation.

\subsection{Approximation results of the stable and unstable DDE limit cycles}\label{Sec_approx_orbits}
 We now illustrate numerically that the GK approximation skills coupled with the accuracy of our effective reduced systems as $\tau$ varies, lead to accurate approximations of  the UPOs to the DDE that unfold through the subcritical Hopf bifurcation, on one hand,  and the stable DDE limit cycles that are produced from the SNO bifurcation, on the other. 
 To prepare for these results, recall that any solution $\T(t)$ to the DDE \eqref{Eq_SS_perturb} emanating from $T_{\textrm{init}}(s)$ ($-\tau \leq s \leq 0$), is approximated by the solution, $\y(t)=(y_1(t),\cdots,y_N(t))^T$, to  the GK system Eq.~\eqref{Eq_Galerkin_SS_6D} that emanates from the projection of  $T_{\textrm{init}}$ onto the Koornwinder polynomials $\mathcal{K}_n^\tau$ (see \eqref{eq:Pn_tilde_prod}),  according to the formula 
\be \label{Eq_SS_reconstruct1}
\T(t)\approx \sum_{j=1}^N  y_j(t),
\ee
see Eq.~\eqref{Eq_App2} and \cite[Sec.~6]{CGLW16}. 

Given a GK system with its dimension, $N$, sufficiently large so that the approximation error in $L^{\infty}$ is small in \eqref{Eq_SS_reconstruct1} (see Table \ref{Table_Error}), we can expect to reach a good approximation of $\T(t)$, at least for $\tau$ sufficiently close to $\tau_c$  (due to Theorem  \ref{Thm_CM_approx}),  by replacing the low-mode amplitudes $y_1(t)$ and $y_2(t)$ by their approximations $x_1(t)$ and $x_2(t)$ from the 2D reduced system \eqref{Eq_EffectiveReduced_ENSO}, and the $y_n(t)$ by $\Psi_{n,\tau}(x_1(t)  \bm{e}_1+x_2(t)  \bm{e}_2)$, for $n \geq 3$.

\begin{figure}[hbtp]
   \centering
\includegraphics[width=0.9\textwidth, height=0.55\textwidth]{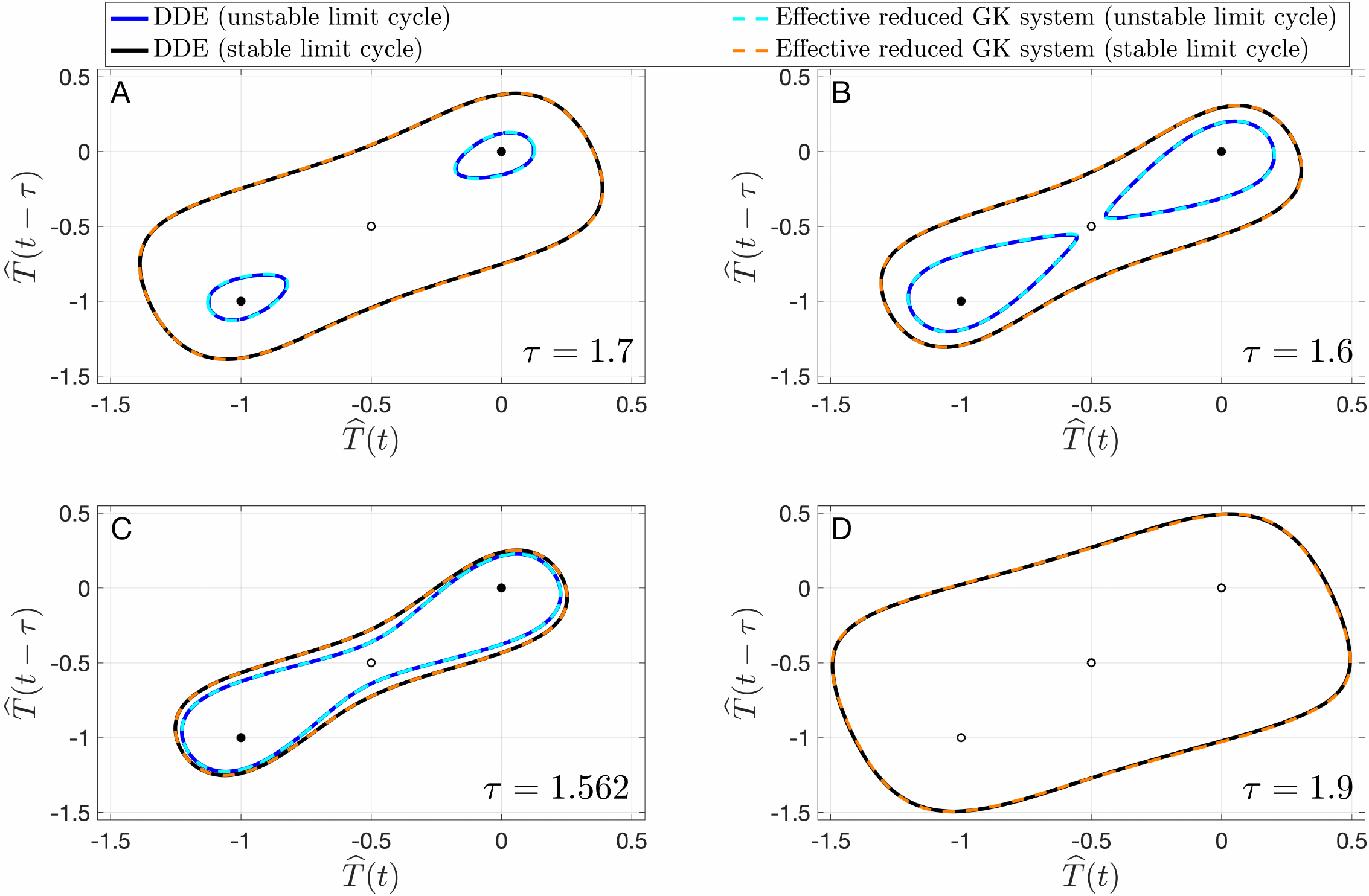}
  \caption{{\bf Approximation results: DDE vs Effective Reduced GK system.}  The solutions are shown in lagged coordinates for $\alpha=0.75$. For different values of $\tau$ as indicated, the solutions to the DDE \eqref{Eq_SS_perturb} are compared to those obtained from the formula of $T^\ast$ given by \eqref{Eq_SS_reconstruct2}, built from the solutions of the 2D reduced GK system Eq.~\eqref{Eq_EffectiveReduced_ENSO} and the parameterization $\Psi_\tau$ given by  \eqref{Eq_Phi}.  The stable (resp.~unstable) DDE limit cycles are shown in black (resp.~blue). The stable (resp.~unstable) limit cycles, obtained from the 2D reduced system after lifting through \eqref{Eq_SS_reconstruct2}, are shown by the  orange (resp.~cyan) dashed curves. The stable (resp.~unstable) equilibria are shown as filled (resp.~empty) circles.  Here, it is sufficient to use a  GK system of dimension $N=6$ to build the parameterization $\Psi_\tau$ and the modes and eigenvalues used in \eqref{Eq_SS_reconstruct2}, in order to achieve such approximation skills.}\label{Fig_Hopf_SS}
\end{figure}

 We consider thus the following approximation formula of $\T(t)$ that is built solely from the solutions of the 2D reduced GK system Eq.~\eqref{Eq_EffectiveReduced_ENSO} and the high-mode parameterization $\Psi_\tau$ given by  \eqref{Eq_Phi}:
\bea \label{Eq_SS_reconstruct2}
&T^*(t)=\sum_{j=1}^N \Big(x_1(t) \boldsymbol{e}_1^j + x_2(t) \boldsymbol{e}_2^j + \sum_{n = 3}^N \Psi_{n,\tau}\big(\x(t)\big) \boldsymbol{e}_n^j \Big), \\
&\mbox{ with } \x(t)=x_1(t)  \bm{e}_1+x_2(t)  \bm{e}_2 \mbox{ solving  Eq.~\eqref{Eq_EffectiveReduced_ENSO}},
\eea
and where $\boldsymbol{e}_{\ell}^j$ denotes the $j^{\mathrm{th}}$ component of the eigenmodes $\boldsymbol{e}_{\ell}$ for $ 1\leq \ell \leq   N$.  Recall that the latter depend on $\tau$ and $N$ as eigenvectors of $A(\tau)$ given by \eqref{eq:A_SS}.

 Figure \ref{Fig_Hopf_SS} shows, for $N=6$, the approximation skills by $T^*(t)$ of periodic orbits to the (perturbed) Suarez and Schopf model \eqref{Eq_SS_perturb}  for four different cases. Three of these cases exhibit  UPOs coexisting with a  stable limit cycle and cover situations for which (i) $\tau$ is close to $\tau_c$ from below ($\tau=1.7$, Fig.~\ref{Fig_Hopf_SS}A), (ii) $\tau$ is close (from below) to $\tau^{\sharp}$ where the homoclinic orbit emerges ($\tau=1.6$, Fig.~\ref{Fig_Hopf_SS}B), and (iii)  $\tau$ is close to $\tau^\ast$ at which the SNO bifurcation occurs ($\tau=1.562$, Fig.~\ref{Fig_Hopf_SS}C). The last case shown corresponds $\tau$ away from $\tau_c$ from above  ($\tau=1.9$, Fig.~\ref{Fig_Hopf_SS}D) for which only a stable limit cycle exists as periodic orbit. 
 In each case, the approximation formula \eqref{Eq_SS_reconstruct2} based on the 2D reduced GK system \eqref{Eq_EffectiveReduced_ENSO} and the parameterization $\Psi_\tau$ enables us to achieve highly accurate approximations of the bifurcating DDE periodic orbits, including the corresponding UPOs and stable limit cycles along with their pronounced nonlinear features, that unfold  from the subcritical Hopf bifurcation and the SNO bifurcation, respectively; see Fig.~\ref{Fig_bif_combo}.

\begin{figure}[hbtp]
   \centering
\includegraphics[width=0.7\textwidth, height=0.35\textwidth]{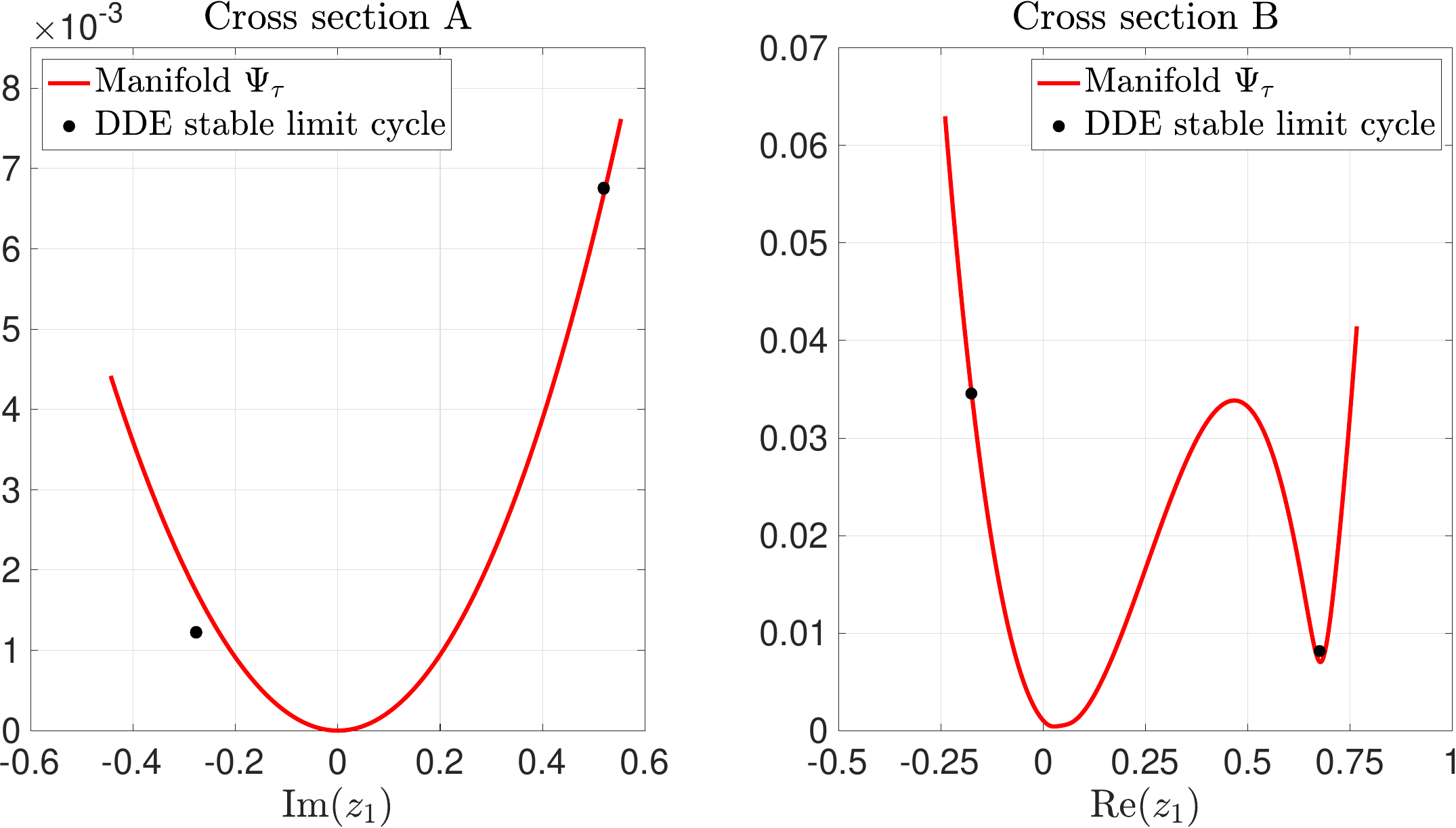}
  \caption{{\bf Visualization of the high-mode parameterization $\Psi_\tau$ used in the effective reduced Eq.~\eqref{Eq_EffectiveReduced_ENSO}.} The curves (in red) are obtained as intersections between the graph of $\varphi$ given by \eqref{Eq_psi} (for $N=6$) with respectively, the vertical plane corresponding to $\mathrm{Re} z_1 \equiv 0$ (cross section A), and that corresponding to $\mathrm{Im} z_1 = a \, \mathrm{Re} z_1 + b$ with $a = -0.63$ and $b = 0.21$ (cross section B). The black dots correspond to the intersections points of the  DDE stable limit cycle from  \eqref{Eq_SS_perturb} (for $\tau=1.9$), with these vertical planes.}   \label{Fig_Phi_visual_SS}
\end{figure}

The reason behind  these impressive skills achieved by the reduced system \eqref{Eq_EffectiveReduced_ENSO}  when $\tau$ varies, 
 lies in a simple property satisfied by the DDE periodic solutions: they sit very close (almost slaved)  to the manifold given as graph of $\Psi_\tau$ built from the GK system \eqref{Eq_Galerkin_SS_6D}, already for $N=6$. This property remains valid even for solution's amplitude (for $\T(t)$) of order one such as shown for $\tau = 1.9$ in  Fig.~\ref{Fig_Phi_visual_SS}, which corresponds to the furthest value $\tau$ to $\tau_c$ examined in  Fig.~\ref{Fig_Hopf_SS}. 

 These comments are made visual in Fig.~\ref{Fig_Phi_visual_SS}. 
This figure shows  the graph of the norm of $\Psi_\tau$,
\be\label{Eq_psi}
\varphi:(z_1, z_2) \mapsto  \|\Psi_\tau(\z) \|_\s, \;\; \z=z_1 \bm{e}_1+z_2 \bm{e}_2,
\ee
as a function of $\mathrm{Re}(z_1)$ and $\mathrm{Im}(z_1)$, where
\be
 \| \w\|_\s=\sqrt{\sum_{j=m_c+1}^N \big| \langle \w,\boldsymbol{e}_j^*\rangle\big|^2}, \mbox{ for any } \w\in H_\s.
\ee

The representation of $\varphi$ as a function of $\mathrm{Re} z_1$ and $\mathrm{Im} z_1$ is made possible since $H_{\c}$ is spanned by 
 eigenvectors that form a complex conjugate pair. The resulting surface, $(\mathrm{Re} (z_1), \mathrm{Im} (z_1)) \mapsto \varphi(\mathrm{Re}(z_1), \mathrm{Im}( z_1))$, is then intercepted by two vertical planes for a better visual inspection of the distance between the stable periodic trajectory from  Eq.~\eqref{Eq_Galerkin_SS_6D} and this surface.  The left and right panels of Fig.~\ref{Fig_Phi_visual_SS} show the curves (in red) obtained as cross sections with, respectively, the vertical plane corresponding to $\mathrm{Re}( z_1) \equiv 0$, and that corresponding to $\mathrm{Im}( z_1) = a \, \mathrm{Re} (z_1) + b$, with $a = -0.63$ and $b = 0.21$.   The intersections between these planes and the stable periodic solution $\Gamma(t)$ to the DDE \eqref{Eq_SS_perturb}
 are shown by the black dots shown in each panel of Fig.~\ref{Fig_Phi_visual_SS}.  
 
It is noteworthy to emphasize that the  good parameterization of the stable periodic solutions to Eq.~\eqref{Eq_Galerkin_SS_6D} are not limited to the choice of cross sections picked up here.  The parameterization skill can be indeed assessed quantitatively  by inspecting e.g.~the parameterization defect  
\bes
 R(t)= \Psi_\tau(\Gamma_1(t), \Gamma_2(t)) -\Gamma_{\s}(t), 
\ees
where $\Gamma_j$  denotes the projection of stable periodic solution $\Gamma(t)$ to the DDE \eqref{Eq_SS_perturb} onto $\boldsymbol{e}_j$ ($j= 1,2$),  while $\Gamma_{\s}(t)$  denotes the projection of $P(t)$ onto $H_\s=\textrm{span}\{\boldsymbol{e}_3,\cdots,\boldsymbol{e}_N\}$.  For the case  shown in Fig.~\ref{Fig_Phi_visual_SS} ($\tau = 1.9$), the time-averaged over one period of the ratio, $\|R(t)\|^2/\|\Gamma_{\s}(t)\|^2$,  is approximately equal  to $2.5\times 10^{-2}$, indicating that $\Gamma(t)$ is almost slaved to the manifold $\Psi_\tau$. 
 This observation explains that the discrepancies observed in  Fig.~\ref{Fig_Phi_visual_SS} between the DDE solution and the 
manifold $\Psi_\tau$ do not affect the approximation skills of the reduced system \eqref{Eq_EffectiveReduced_ENSO} shown in Fig.~\ref{Fig_Hopf_SS}. Indeed, the fraction of the energy contained in the parameterized modes is so small compared to that contained in the resolved ones, that such discrepancies do not impact the reduced system's approximation skills.

\subsection{Model error estimates}\label{Sec_model_error}
To complement the numerical results of the previous section, 
we present now some basic error estimates for the effective reduced GK system \eqref{Eq_GK_ENSO}. 
As shown in Proposition \ref{Prop_error} below, the parameterization defect associated with the manifold $\Psi_\tau$, is the main controlling factor. Although we present the result in the context of the Suarez-Schopf model, the same type of calculations can be performed for  GK approximations of DDEs with polynomial nonlinearities. The result is also not limited to reduced systems constructed from the center-unstable manifold parameterizations such as given by \eqref{Eq_phi_leadingorder} or \eqref{Eq_Phi}. For this reason, we present the calculation for reduced system based on an arbitrary parameterization function $\varphi_\tau\colon H_{\c} \rightarrow H_{\s}$.  

Let us first recall that the $N$-dimensional GK system \eqref{Eq_Galerkin_SS} of the Suarez-Schopf DDE \eqref{Eq_SS_perturb} is given by 
\be \label{Eq_GK_recall}
\frac{\d \y}{\d t} = A (\tau)\y + G_2(\y) + G_3(\y),
\ee
with the matrix $A(\tau)$ and the nonlinear terms given by \eqref{eq:A_SS} and \eqref{Eq_G_decomp_SS}, respectively. As in Sec.~\ref{Sec_GK_centerman}, we assume that $A(\tau)$ is diagonalizable over $\mathbb{C}$, which is observed to be true for all the numerical results reported in this study for this DDE model. Recall that $H_{\c}$ and $H_{\s}$ are the eigen subspaces spanned by the low modes and the high modes, respectively. 

We will also make use of the following notations. Given a solution $\y(t)$ of \eqref{Eq_GK_recall} evolving in $\mathbb{R}^N$, denote by $\y_\c=\Pi_\c \y = \sum_{j=1}^{m_c} z_j \boldsymbol{e}_j$ the low-mode projection of $\y$ and by $\y_\s=\Pi_\s \y = \sum_{j=m_c+1}^{N} z_j \boldsymbol{e}_j$ the high-mode projection of $\y$. Note that the vector $\z = (z_1, \ldots, z_N)^T$ is nothing else than the vector $\y$ under the eigenbasis of $A(\tau)$. Let $P = [\boldsymbol{e}_1, \ldots, \boldsymbol{e}_N]$ be the $N\times N$ matrix whose columns consist of the eigenvectors of $A(\tau)$. We have $\y = P\z$. 

We introduce also $P_{\c}$ as the matrix made of the  first $m_c$ columns of $P$, namely $P_\c=  [\boldsymbol{e}_1, \ldots, \boldsymbol{e}_{m_c}]$. Similarly, $P_{\s } = [\boldsymbol{e}_{m_c+1}, \ldots, \boldsymbol{e}_N]$, $\z_{\c} = (z_1, \ldots, z_{m_c})^T$ and $\z_{\s} = (z_{m_c + 1}, \ldots, z_{N})^T$. We have then $\y_{\c} = P_{\c} \z_{\c}$ and $\y_{\s} = P_{\s} \z_{\s}$. Finally, let $Q$ be the $N\times N$ matrix whose columns consist of the eigenvectors of $A^*$. 
We have then 
\be
Q^* P = I_{N\times N},
\ee
where $Q^*$ denotes the complex conjugate transpose. Similarly as for $P_\c$, the matrix $Q_\c$ denotes the matrix made of the  first $m_c$ columns of $Q$.

\bprop\label{Prop_error}
Consider $\y(t)$ to be a solution to Eq.~\eqref{Eq_GK_recall} over a given interval $[0,t_m]$. Assume that $A(\tau)$ is diagonalizable over $\mathbb{C}$. Let $\varphi_\tau\colon H_{\c} \rightarrow H_{\s}$ be a parameterization of $\y_\s$ in terms of $\y_\c$. Assume that $\y(t)$ satisfies 
\be\label{Assumption_R}
\|\y(t)\| \le R,  \qquad \|\y_\c(t) + \varphi_\tau(\y_{\c} (t))\| \le R, \qquad t\in[0, t_m].
\ee
Then there exists a constant $C(R)>0$ such that
\be \label{Eq_error_est}
\overline{\left\|\dot{\z}_\c-  \Lambda_{\c} \z_{\c} -Q_\c^* G_2\Big( P_{\c} \z_{\c} + \varphi_\tau(P_{\c}\z_{\c} ))\Big)  \right\|^2 } \leq C(R) \overline{ \|P_{\s} \z_{\s} - \varphi_\tau(P_{\c}\z_{\c}) \|^2},
\ee
where $\overline{(\cdot)}$ denotes the time average over $(0,t_m)$ and $\Lambda_{\c}$ denotes the $m_c \times m_c$ diagonal matrix made of the eigenvalues $\lambda_1, \ldots, \lambda_{m_c}$ of $A(\tau)$.
\eprop

{\it Proof.} The desired result \eqref{Eq_error_est} follows from a direct calculation by comparing the vector field of the reduced system associated with $\varphi_\tau$ (see \eqref{Eq_reduced_w} below) and the projection of the GK vector field in \eqref{Eq_GK_recall} onto $H_{\c}$. Indeed, note that the GK system \eqref{Eq_GK_recall} can be rewritten as:
\be
\dot{\z}= \Lambda(\tau) \z+ Q^* \left(G_2( P \z ) + G_3( P \z )\right),
\ee
where $\Lambda (\tau)= Q^* A(\tau) P$.

By projecting the above equation onto the low modes, we have
\be
\dot{\z_{\c} }= \Lambda_{\c} (\tau) \z_{\c} + Q_{\c}^* \left( G_2( P \z)  + G_3( P \z )\right).
\ee
On the other hand, the reduced system associated with $\varphi_\tau$ is
\be \label{Eq_reduced_w}
\dot{\w}= \Lambda_{\c} (\tau) \w +  Q_\c^* \big( G_2( P_{\c}\w + \varphi_\tau(P_{\c} \w)) + G_3( P_{\c}\w + \varphi_\tau(P_{\c}\w)) \big).
\ee
We define the model error as
\bea \label{Eq_model_err1}
& R_T(\z_\c, \varphi_\tau)=\int_{0}^{t_m} \left\|\frac{\d \z_{\c} }{\d t} -  \Lambda_{\c} \z_{\c} -Q_\c^* \Big( G_2 \Big( P_{\c}\z_{\c} + \varphi_\tau(P_{\c}\z_{\c})\Big) + G_3\Big( P_{\c}\z_{\c} + \varphi_\tau(P_{\c}\z_{\c}) \Big) \Big) \right\|^2 \d t \\
& = \int_{0}^{t_m} \left\| Q_\c^* \Big( G_2(\y) + G_3(\y) - G_2 \Big( P_{\c}\z_{\c} + \varphi_\tau(P_{\c}\z_{\c})\Big) - G_3 \Big( P_{\c}\z_{\c} + \varphi_\tau(P_{\c}\z_{\c}) \Big)  \Big) \right\|^2 \d t \\
& \le C \int_{0}^{t_m} \left\| G_2(\y) + G_3(\y) - G_2\Big( P_{\c}\z_{\c} + \varphi_\tau(P_{\c}\z_{\c}) \Big) - G_3 \Big( P_{\c}\z_{\c} + \varphi_\tau(P_{\c}\z_{\c}) \Big)  \right\|^2 \d t.
\eea
Let us introduce the notation 
\be
\w_\s= \varphi_\tau(P_{\c}\z_{\c}),
\ee
and recall also that $\y_{\c}= P_{\c} \z_{\c}$ and $\y_\s= P_{\s} \z_\s$. Then, by the definitions of $G_2$ and $G_3$ given in \eqref{Eq_G_decomp_SS}, we have 
\bea\label{G2_estim}
G_2(\y) - G_2\Big( P_{\c} \z_{\c} + \varphi_\tau(P_{\c}\z_{\c}) \Big) &= - 3 T_{+}  \Big( \Big( \sum_{n=1}^N y_{\c,n} + \sum_{n=1}^N y_{\s,n}\Big)^2 - \Big(\sum_{n=1}^N y_{\c,n} + \sum_{n=1}^N w_{\s,n}\Big)^2 \Big) \bm{\nu}_N,\\
&= - 3 T_{+} \Big(\sum_{n=1}^N(2 y_{\c,n} + y_{\s,n} + w_{\s,n}) \Big)  \Big(\sum_{n=1}^N \big(y_{\s,n} - w_{\s,n}\big) \Big)  \bm{\nu}_N,
\eea
and
\bea\label{G3_estim}
G_3(\y) - G_3\Big( P_{\c} \z_{\c} + \varphi_\tau(P_{\c}\z_{\c}) \Big) &= - \Big( \Big( \sum_{n=1}^N y_{\c,n} + \sum_{n=1}^N y_{\s,n}\Big)^3 - \Big(\sum_{n=1}^N y_{\c,n} + \sum_{n=1}^N w_{\s,n}\Big)^3 \Big) \bm{\nu}_N \\
&= - \mathcal{N}(\y,\w_\s) \Big(\sum_{n=1}^N \big(y_{\s,n} - w_{\s,n}\big) \Big)  \bm{\nu}_N,
\eea
where
\beas
\mathcal{N}(\y,\w_\s) &= 3 \Big(\sum_{n=1}^N y_{\c,n}\Big)^2 + \Big(\sum_{n=1}^N y_{\s,n}\Big)^2 + \Big(\sum_{n=1}^N w_{\s,n}\Big)^2 \\
& \qquad + 3 \Big(\sum_{n=1}^N y_{\c,n}\Big) \Big(\sum_{n=1}^N (y_{\s,n} + w_{\s,n})\Big) + \Big(\sum_{n=1}^N y_{\s,n}\Big)\Big(\sum_{n=1}^N w_{\s,n}\Big).
\eeas
Finally, by using \eqref{G2_estim} and \eqref{G3_estim} in \eqref{Eq_model_err1}, we conclude, in virtue of assumption \eqref{Assumption_R}, that 
\bes
R_T(\z_\c, \varphi_\tau) \leq C(R) \int_{0}^{t_m} \left| \sum_{n=1}^N \big(y_{\s,n} - w_{\s,n}\big) \right |^2\d t   \leq C(R) \int_{0}^{t_m} \|\y_{\s}(t) - \w_{\s} (t) \|^2\d t.\hspace{2.5cm} \square
\ees


\section{Tipping Solution Paths and ENSO variability}\label{Sec_ENSO_var}

\subsection{Transition paths and nonlinear building blocks of temporal variability}\label{Sec_nln_blocks}

The plausibility of a  stochastic forcing as a mechanism for ENSO irregularity has been argued in many studies. 
Physical origins of such a forcing include large-scale synoptic atmospheric transients such as the Madden Julian Oscillation  \cite{batstone2005characteristics} or westerly wind bursts \cite{fedorov2002response}. Dynamically, the idea is to explicitly separate the slow and fast modes in the atmosphere and add the latter to simple deterministic models as a stochastic forcing term.
Other sources of stochasticity include processes associated with atmospheric/moist convective disturbances whose timescales can vary from hours to weeks.  
Typically, the additional fast-mode random forcing disrupts the slow scales and convert the original periodic or damped oscillation supported by the deterministic model into an irregular one. Such mechanisms have been previously advocated through a combination of observational data and a hierarchy of ENSO models including data-driven models \cite{penland1995optimal,kondrashov2005hierarchy,CKG11,chen2016diversity};  PDE models \cite{blanke1997estimating,chen2018observations,eckert1997predictability,thual2016simple,roulston2000response,zavala2003response}, and conceptual models \cite{chen2022multiscale,chen2023rigorous}.

Here, we propose a novel scenario for the fabric of ENSO variability. Our approach is based on  (i) the dynamical insights gained from the bifurcation diagram and phase portrait of the deterministic dynamics as summarised by Fig.~\ref{Fig_bif_combo}, and (ii) the design of stochastic model's disturbances interacting with  the underlying nonlinear invariant sets (UPOs, stable limit cycles, and homoclinic orbit) and the topology  they form across a  $\tau$-interval $[\tau_0,\tau_1]$ containing $[\tau^\ast, \tau^{\sharp}]$, where $\tau^\ast$ (resp.~$\tau^{\sharp}$)  corresponds to the SNO  (resp.~homoclinic) bifurcation point.   Since these invariant sets are $\tau$-dependent, as dictated by the bifurcation diagram of Fig.~\ref{Fig_bif_combo}, we naturally consider these rapidly varying disturbances to be superimposed to slowly varying lagged effects. Physically, such slow variations can be envisioned as resulting from slight variations  occurring in key properties of the equatorially trapped oceanic waves such as their transit time of the Pacific ocean. 
 
 In that respect, we propose the following stochastic model that we write in the perturbed variable $\theta(t)=T - T_{+}$:
\begin{subequations} \label{Eq_stoch_model}
\begin{align}
&\d \theta (t)= \Big(a\theta (t)- \alpha \theta(t - \tau(t)) - b\theta^2(t) - \theta^3(t)\Big) \d t + \frac{\sigma}{1+|\theta^2(t)|}\d W_t, \label{Eq_stoch_model_a}\\
&\mbox{with }\tau(t)=\tau_0+\epsilon t,  \qquad \tau(t)\in [\tau_0,\tau_1],\label{Eq_taut} 
\end{align}
\end{subequations} 
where $a=1 - 3 T_{+}^2$, $b=3 T_{+} $, and $\epsilon, \sigma\geq 0$ while $W_t$ denotes a Brownian motion.
The noise term is a Lorentzian function favouring larger noise when $\theta\approx 0$ while damping the noise effects when $\theta$ gets too large. It is interpreted in the sense of It\^{o} to fix ideas \cite[Chapter 4]{gardinerbook}.  Due to this noise term, a stochastic path solving \eqref{Eq_stoch_model} is more likely to meander away from $T_+$ than around it, as time flows. Since $T_+$ corresponds for the Suarez and Schopf model to the basic steady state associated with an El Ni\~no event, using such a nonlinear noise favours, intuitively, a less frequent occurrence  of  metastable El Ni\~no events corresponding to visiting a neighborhood of $T_+$ than if  Eq.~\eqref{Eq_stoch_model} would be driven by a pure white noise.\footnote{The nonlinear noise favours also a less frequent occurrence of El Ni\~no events with unrealistic long duration visits that would consist of e.g.~meandering near $T_+$ for a too long amount of time.} As discussed below, Eq.~\eqref{Eq_stoch_model} supports actually the occurrence of other types of El Ni\~no events that do not correspond to a metastable visit of $T_+$ but rather to experiencing excursions across the bifurcating solutions of the deterministic equation as identified in Sec.~\ref{Sec_Bif_SNO}.

To show evidence of these other types of El Ni\~no events, we performed an integration of Eq.~\eqref{Eq_stoch_model} for $\alpha=0.75$, $\sigma=0.2$, with  $\tau$ evolving according to Eq.~\eqref{Eq_taut} over an interval $[0,{t_m}]$, with $\tau_0=1.45$, $\epsilon=8.4\times 10^{-4}$ , and ${t_m}=237.8$ in time unit of the model. With these parameters,  we have $\tau(t_m)=\tau_1=1.65$, and thus  the interval $[\tau_0,\tau_1]$, within which $\tau(t)$ varies, contains both the SNO critical value $\tau^\ast$ and the homoclinic critical value $\tau^{\sharp}$; see caption of Fig.~\ref{Fig_bif_combo} for numerical values of $\tau^\ast$ and $ \tau^{\sharp}$. The rationale behind this numerical setup is guided by a simple intuition. Indeed, as the delay parameter $\tau$ slowly drifts across a portion of the bifurcation diagram shown in Fig.~\ref{Fig_bif_combo} and when the noise strength parameter $\sigma$ is not too small, it is expected that in the course of integration of Eq.~\eqref{Eq_stoch_model}, its solution would ``hop'' through the different branches of the bifurcation diagram and thus reveal fingerprints of the various invariant sets (UPOs, stable limit cycles, steady states) of the deterministic flow. We call such a solution a Tipping Solution Path (TSP).

\begin{figure}[tbph]
   \centering
\includegraphics[width=0.95\textwidth, height=0.4\textwidth]{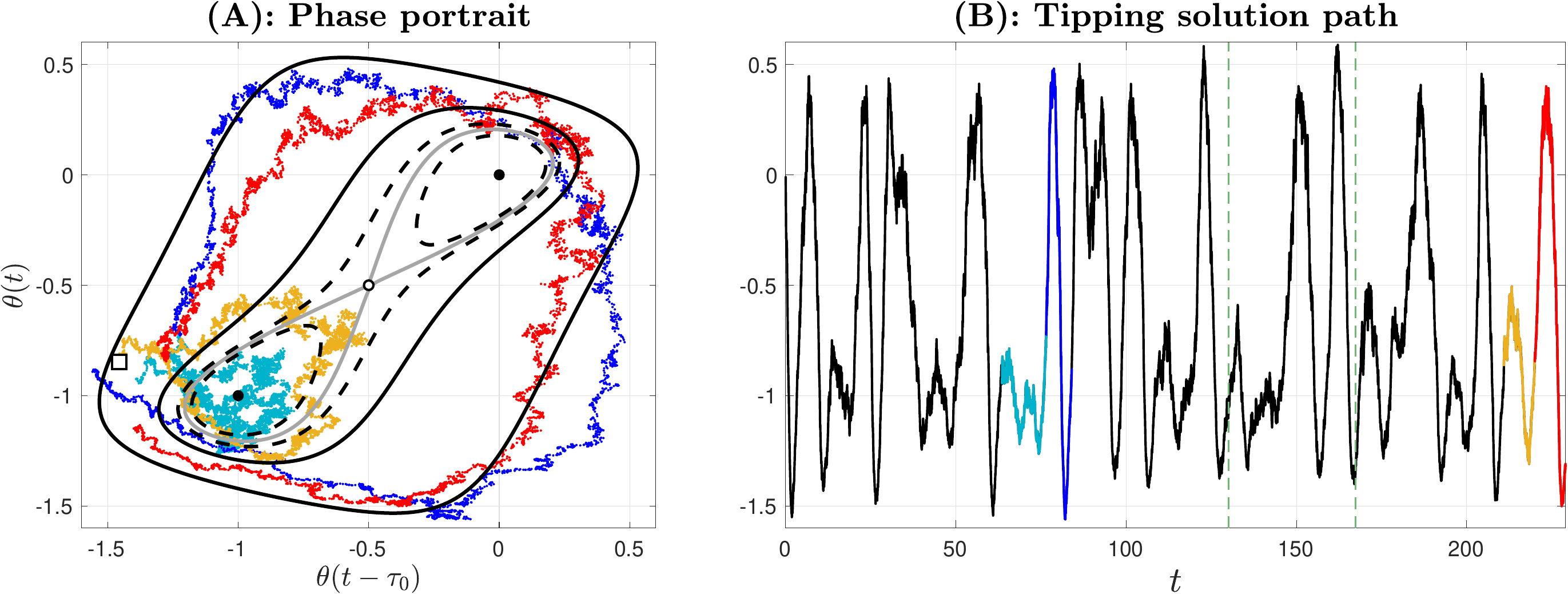}
  \caption{{\bf Tipping solution path and phase portrait visualization of two of its transition snippets.} One transition  snippet emanates near the La Ni\~na steady state $T_{-}$ (of negative value) and spends some time meandering around it (turquoise path in panel A) after escaping it via a nearly cyclic noisy trajectory (blue path in panel A). The corresponding temporal pattern is shown on $\theta(t)$ in panel B with the same color coding.  The other transition snippet emanates from a location marked by the empty square in panel A and meanders around the UPO surrounding $T_{-}$ before escaping via a nearly cyclic noisy trajectory (red path in panel A). The corresponding episode is shown in the time domain with the same color coding in panel B. The two vertical green dashed  lines in panel B mark the time-instants at which $\tau(t)=\tau^*$ ($t = 130$) and $\tau(t)=\tau^{\sharp}$ ($t = 167.38$). In panel A, besides the stochastic trajectory, are shown several invariant sets (UPOs, stable limit cycles, and homoclinic orbit) for the deterministic model, for different $\tau$-values as pointed out in the text.}\label{Fig_4segments}
\end{figure}

\begin{figure}[tbph]
   \centering
\includegraphics[width=0.95\textwidth, height=0.4\textwidth]{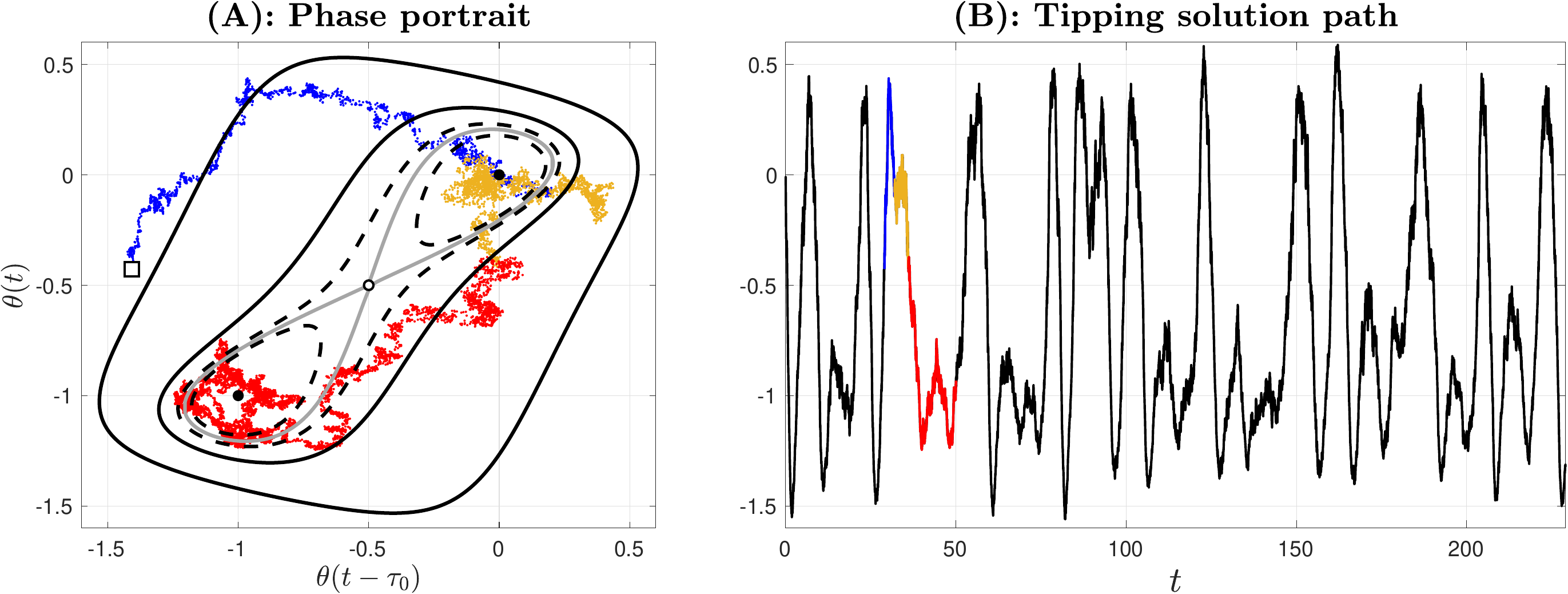}
  \caption{{\bf A transition snippet going through a neighborhood of $T_{+}$ before landing near $T_{-}$}. The empty square in panel A marks the beginning of the blue path. See Text.}\label{Fig_hopping}
\end{figure}

In the literature, different tipping mechanisms have been identified that include bifurcation-induced tipping, noise-induced tipping, and rate-induced tipping \cite{ashwin2012tipping,feudel2018multistability,kuehn2011mathematical}. It is known that under different setups of the external forcing and/or parameter drifting scenarios, a given model may experience different types of tipping; see e.g.~\cite[Section 4]{ashwin2012tipping}.  A mathematical characterization of the different types of  tipping phenomenon in terms of the model parameters in  Eq.~\eqref{Eq_stoch_model} and drifting scenarios, requires a separate study. 
Instead, our focus is on the impact on the solution's variability of the presence of  tipping points involving UPOs and in that regard, the parameter setting and drifting scenario \eqref{Eq_taut} chosen above have been found to be physically insightful.
In particular, as explained below (see Figs.~\ref{Fig_4segments} and \ref{Fig_hopping}), we observed that even during  time intervals preceding the crossing of the SNO critical value $\tau^*$, the solution's variability is affected by invariant sets such as limit cycles that emerge only for $\tau>\tau^*$.

We numerically observed that this ability of the stochastic model to ``feel'' the deterministic model's dynamics at a given $\tau$-value  prior $\tau(t)$ crosses this value, is actually  independent of how the noise is interpreted (It\^{o} vs.~Stratonovich).
 \footnote{In the Stratonovich sense (\cite[Section 4.3.6]{gardinerbook}), an additional term $\sigma^2 \theta(t)(1+\theta^2(t))^{-3} \d t$ is added to Eq.~\eqref{Eq_stoch_model_a}. In the present parameter setting, we have checked that this extra term has little influence on the dynamics  as it scales here with $\sigma^2$.} 
Thus, noise can e.g.~excites parameter regimes with oscillatory dynamics as $\tau(t)$ approaches $\tau$-values supporting oscillations for the unperturbed model. In terms of solution's variability, this can be reflected on trajectory segments of the stochastic solution exhibiting then ``fingerprints'' of these underlying nonlinear oscillations.

Thus, noise-sustained oscillations associated with the SNO critical transition play an important role in the fabric of the temporal variability of a TSP solving Eq.~\eqref{Eq_stoch_model}. Their role is not exclusive though as other nonlinear invariant sets (steady states, homoclinic orbits) have also their share in the  TSP's variability.
In that respect, we discuss below how the whole set of these nonlinear building blocks allow us to break down ENSO-like events from the model \eqref{Eq_stoch_model}.

Such a TSP is shown in Fig.~\ref{Fig_4segments}B. To understand the dynamical origins behind the time-variability displayed by this TSP, we focus on  a few episodes (snippets) as marked in color in Fig.~\ref{Fig_4segments}B  that we map, with the same color coding, to the phase portrait shown in Fig.~\ref{Fig_4segments}A. 
To help interpretation, these snippets are superimposed on a few invariant curves of the deterministic model \eqref{Eq_SS_perturb}\footnote{as calculated from our 2D reduced equation; see  Secns.~\ref{Sec_Bif_SNO} and \ref{Sec_approx_orbits}.} encountered as $\tau$ evolves from  $\tau_0$  to $\tau_1$ according to \eqref{Eq_taut}: a UPO at $\tau=1.5607$ encompassing the homoclinic orbit, two UPOs symmetric  with respect to the saddle point (for $\tau=\tau_1=1.65$), and  two stable limit cycles,   shown as two closed black (plain) curves
whose the one with smaller diameter corresponds to  $\tau=1.65$ while the other corresponds to $\tau=2$.

The snippets marked by the turquoise and blue consecutive segments of the TSP as shown in Fig.~\ref{Fig_4segments}B, correspond to a quiet episode (turquoise) followed by a large excursion (blue) reminiscent with what has been observed for certain El Ni\~no events as recorded by the Ni\~no 3.4 SST index. Such  a large excursion following a certain stillness evokes the major El Ni\~no event that  occurred in Dec-Jan 2016 \cite{BAMS_2017observing}. It is worthwhile noting that the Ni\~no 3.4 SST index during the four years preceding this event exhibited indeed temporal patterns resembling the episode marked in turquoise; cf.~e.g.~the global reanalysis product of Ni\~no 3.4 index as documented by the E.U.~Copernicus Marine Service Information \cite{EU_Nino34}.
The phase portrait reveals  that these consecutive events correspond,  in our model, to a meandering of the solution path near what corresponds to a La Ni\~na steady state $T_{-}$ (turquoise path in Fig.~\ref{Fig_4segments}A)\footnote{To avoid multiple notations, we still denote by $T_{-}$ the steady state to Eq.~\eqref{Eq_SS_perturb} that corresponds to $T_{-}$ given by \eqref{Eq_steady_states}.}  followed by a nearly cyclic orbit (blue path in  Fig.~\ref{Fig_4segments}A).

The snippets shown  in yellow and red in Fig.~\ref{Fig_4segments}B emphasize a different mechanism   leading to an El Ni\~no event. Indeed, by inspecting the corresponding path in Fig.~\ref{Fig_4segments}A, we observe that while the  El Ni\~no excursion corresponds still to a nearly cyclic orbit  (red), its ignition is now preceded by a path's wandering around an UPO surrounding  itself 
the La Ni\~na steady state. 
Due to their  behavior starting either from the La Ni\~na steady state or an UPO surrounding it, before engaging into a quasi cyclic excursion, we call such episodes as Transition Snippets.   

Such transition snippets provide the generic building blocks at work for the production of El Ni\~no-like patterns displayed by the solution to Eq.~\eqref{Eq_stoch_model} such as shown in Fig.~\ref{Fig_4segments}B.  Our analysis reveals thus that a solution path wandering around  an UPO  surrounding the La Ni\~na steady state or having a closer meandering around it can be interpreted as an early warning signal \cite{scheffer2009early} for an El Ni\~no-like pattern to occur in this model.  

 Other episodes displayed by the TSP's temporal patterns correspond to a path connecting the El Ni\~no with the  La Ni\~na equilibria. 
Such an event is shown in  Fig.~\ref{Fig_hopping} with color coding marking the following three stages: (i) path escape 
 towards the El Ni\~no equilibrium (blue path), (ii) path meandering in a neighborhood of that equilibrium (yellow path), and (iii) path escape from that neighborhood towards  the La Ni\~na equilibrium (red path). 

The mechanisms of fabric of the TSP's temporal patterns reported here  are therefore diverse, 
involving various scenarios such as noise-sustained oscillations associated with the SNO critical transition and its constitutive UPOs and stable limit cycles, or  transitions between steady states. Due to the presence of (nonlinear) periodic orbits these mechanisms are characteristic of more intricate interactions between noise and nonlinear effects than classically encountered in  noise-induced tipping dynamics between multiple equilibria; see e.g.~\cite{gammaitoni1998stochastic} and \cite[Chapter 6]{Horsthemke_book}. 
In particular, the good capture of such interactions by stochastic reduced models remains to be explored. We mention that stochastic invariant manifolds techniques  \cite{CLW15_vol2,chekroun2023transitions} and their extension \cite{CLM23} should be useful in that respect.

 Finally, it worth mentioning that such mechanisms are not dependent on the noise path used to drive Eq.~\eqref{Eq_stoch_model_a} but rather conditioned to the choice of $\tau_0$ in  Eq.~\eqref{Eq_taut} and $\sigma$ in Eq.~\eqref{Eq_stoch_model_a}. Choosing $\tau_0$  too far below $\tau^*$ and/or $\sigma$ too small could lead to completely different solution's temporal behavior with for instance more frequent transitions between equilibria and very rare excitation of oscillatory behavior. 
While allowing $\tau$ in Eq.~\eqref{Eq_stoch_model} to drift is motivated by physical considerations as mentioned above, it is  not excluded though, for well-calibrated model's parameters, that stochastic solutions for frozen $\tau$-value could share similar temporal variability attributes than those reported here.  Regardless, drifting or not, the nonlinear invariant sets (UPOs, stable limit cycles, steady states) of the deterministic flow play a central role in shaping the temporal patterns exhibited by the TSP.  We have been though voluntarily evasive on the role of the homoclinic orbit in the  fabric of the TSP variability. The next section clarifies this point.

\subsection{Decadal variability and homoclinic orbit}\label{Sec_homocline}

 To decipher the role of the homoclinic orbit in the fabric of the TSP variability, we need to quantify this variability in more precise terms.
To do so, we analyze the spectral content of a TSP  produced over a long-time integration of Eq.~\eqref{Eq_stoch_model} consisting of $10^6$ time steps with $\delta t=2 \times 10^{-3}$, in which $\tau(t)$ is now allowed to oscillate in a periodic fashion between $\tau_0$ and $\tau_1$ in \eqref{Eq_taut}. The values of these and other model's parameters are the same as in Sec.~\ref{Sec_nln_blocks}. 
To make a physical sense of the variability displayed by the TSP, we adopt the following conversion formula \cite{boutle2007nino}  $s=  t \Delta/\tau$ with $\Delta=349,$ and $\tau=1.7$, in which $s$ denotes the physical time in year, while $t$ denotes the non-dimensional time unit of the DDE model. To analyse the frequency content, we use 
 the data-adaptive harmonic decomposition (DAHD) method \cite{chekroun2017data} that is a signal processing method which has been employed in many applications to analyze the temporal variability and extract coherent patterns of complex (multivariate) time series; see e.g.~\cite{MASIE_paper,KCYG17,kondrashov2018data,KCB18}.
 Said in simple terms, the method allows for performing a Fourier-like analysis 
 from time-lagged correlations while extracting empirical modes of variability that are naturally ranked per Fourier frequency \cite[Theorem V.1]{chekroun2017data}. 
A simple projection of the signal onto a  selected group of such modes over a given frequency range allows in turn for calculating in the time domain the contribution of that range within the original signal.

\begin{figure}[hbtp]
   \centering
\includegraphics[width=0.95\textwidth, height=0.3\textwidth]{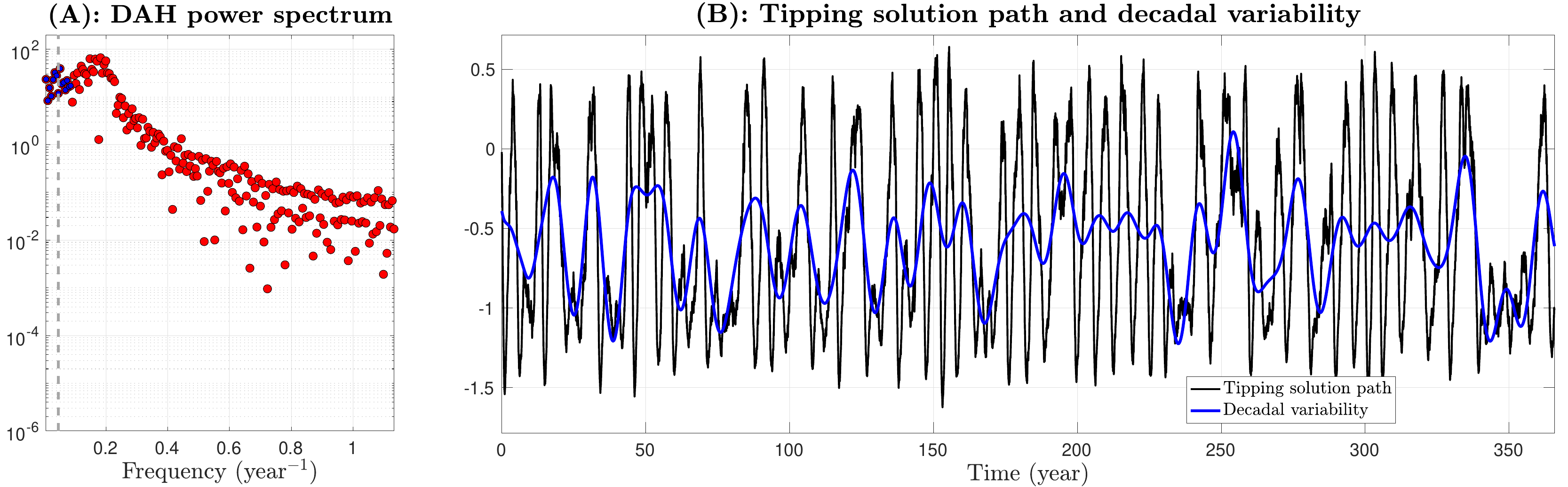}
  \caption{{\bf DAH power spectrum and decadal variability}. The decadal peak located at $\approx$ 21 yr is marked by the vertical grey dashed line in panel A. The ENSO peak is located around $\approx 5.5$ yr. The projection of the TSP onto the subspace spanned by the DAH modes (DAHMs) associated with the decadal variability (blue dots in panel A) is shown by the blue curve in panel B.}\label{Fig_decadal}
\end{figure}

 We applied DAHD to the aforementioned TSP simulated over a time window of approximately $1,100$ yrs. The DAH power spectrum reveals two distinct broadband peaks:  a $21$-yr peak associated with decadal variability, and a $5.5$-yr peak associated with  ENSO interannual variability, corresponding in Fig.~\ref{Fig_decadal}A, to the small bump made of blue dots for decadal variability, and to the more energetic one located to its right, for the interannual one. The projection of the TSP onto the subspace spanned by the DAH modes (DAHMs) associated with the decadal variability (blue dots in Fig.~\ref{Fig_decadal}A) is shown by the blue curve in Fig.~\ref{Fig_decadal}B.

\begin{figure}[hbtp]
   \centering
\includegraphics[width=0.95\textwidth, height=0.3\textwidth]{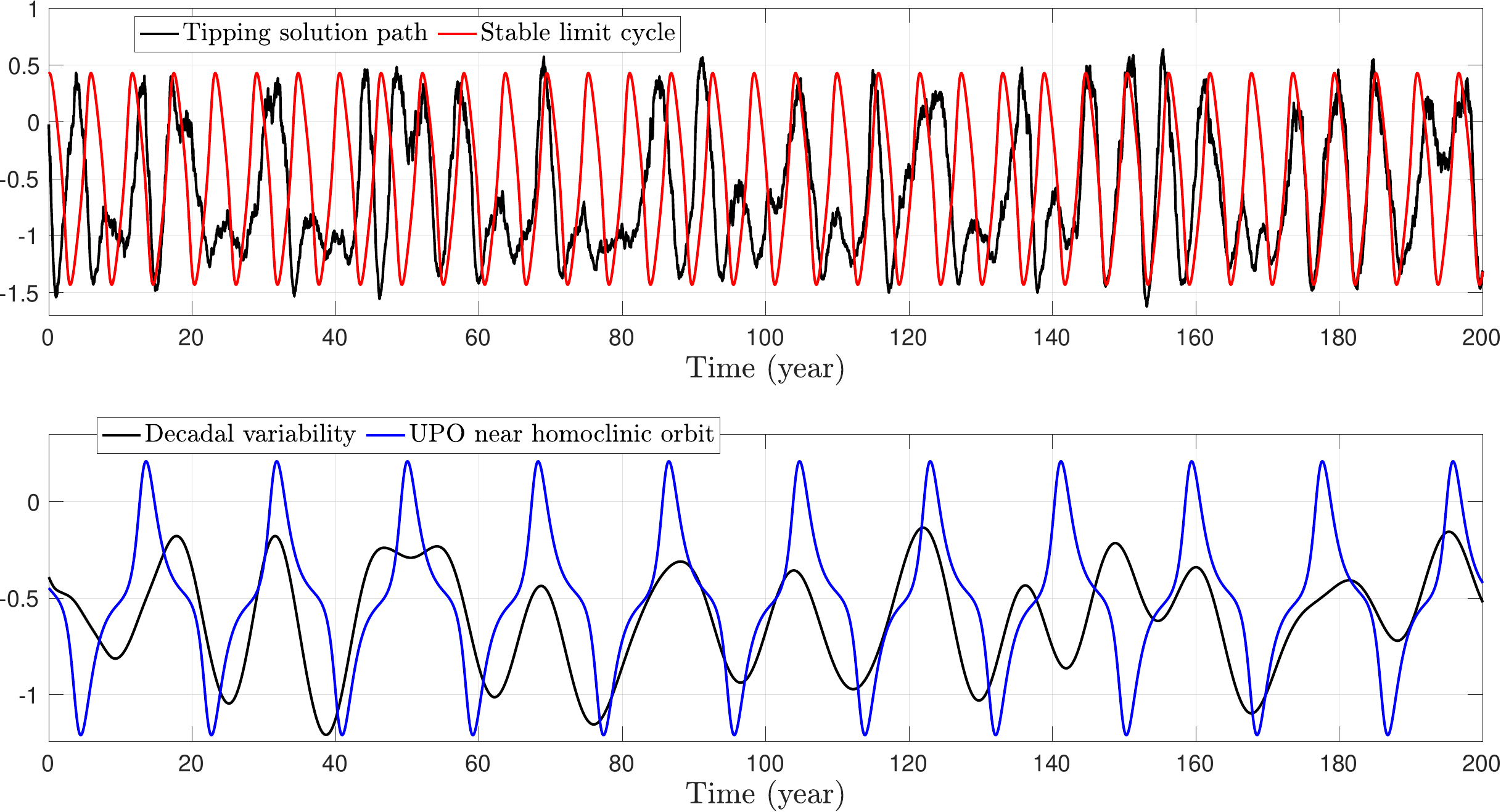}
  \caption{{\bf Interannual  and Decadal variability}. The black curve shown in the lower panel corresponds to the decadal variability time series shown in  blue in the right panel of Fig.~\ref{Fig_decadal}. The blue curve shown here corresponds to a UPO of period $18.23$-yr obtained for $\tau = 1.5827$, which  lies close to the homoclinic orbit. The upper panel shows the TSP in black and a stable limit cycle in red of period $5.78$-yr obtained for $\tau = 1.7689$.}\label{Fig_decadal2}
\end{figure}

Decadal variability of ENSO has been documented in recent studies \cite{dieppois2021enso}.
The origin of this decadal variability in our model \eqref{Eq_stoch_model} can be traced back to UPOs that are located close to the homoclinic orbit, i.e.~for $\tau$ close to $\tau^{\sharp}$ from below; see blue curves in inset B of Fig.~\ref{Fig_bif_combo}. 
Recall that the 2D reduced GK system \eqref{Eq_EffectiveReduced_ENSO} allows for a simple computation of the UPOs of the DDE \eqref{Eq_SS_perturb} by a simple backward integration of Eq.~\eqref{Eq_EffectiveReduced_ENSO} followed by a lifting of the orbit via the formula \eqref{Eq_SS_reconstruct2}. Using this approach, UPOs closer to the homoclinic orbit such as shown in Fig.~\ref{Fig_bif_combo}B can be easily probed. Among these UPOs, we found that the UPO with periodicity of $18.23$-yr (blue curve in Fig.~\ref{Fig_decadal2}-lower panel,  corresponding to $\tau = 1.5827$ in the deterministic model) provides a good synchronization with the temporal (modulated) patterns exhibited by the decadal variability extracted from the DAHD (black curve in Fig.~\ref{Fig_decadal2}-lower panel).  Thus, one can infer that the homoclinic bifurcation occurring at $\tau=\tau^{\sharp}$ is responsible for the emergence of UPOs that provide the backbone of the decadal variability exhibited by our TSP. As shown in the top panel of Fig.~\ref{Fig_decadal2}, the stable limit cycles are responsible for the more regular oscillatory events occurring on the TSP on an interannual timescale.

\section{Concluding Remarks} \label{Sec_conclusion}

Thus, our approach based on GK approximations and higher-order approximations to center-unstable manifolds to analyze bifurcations in the Suarez and Schopf model allowed us to exhibit the nonlinear building blocks (UPOs, stable limit cycles, homoclinic orbits, steady states) at play in the fabric of a rich temporal variability exhibited by solution paths drifting through the underlying  bifurcations in presence of noise. The drift corresponds to slow changes in e.g.~the traveling time of equatorially trapped waves (i.e.~slow drift in $\tau$ through its critical values) while the noise term in \eqref{Eq_stoch_model} can be interpreted as a nonlinear coupling term between ocean and atmosphere.   

This new mechanism identified in this study to produce interannual and decadal variability of ENSO-like patterns deserves more explorations, in particular in  ENSO models  including more physics while still tractable mathematically, such as the Cane-Zebiak models and the like; see \cite{Jin_al93_part1,Jin_al93_part2,cao2019mathematical,Tantet_al_ENSO}.
In that respect, exploring the effects of distributed delays in the Suarez and Schopf could be insightful as suggested in \cite{falkena2019derivation}. The approach presented in this work applies to such models where the integral terms of the form $\int_{t-\tau}^t x(s) \d s$ are simply approximated by $\tau y_0(t) - \tau \sum_{n=1}^{N-1} y_n(t)$, in the corresponding GK approximations; see \eqref{eq:G}. Similarly the inclusion of multiple delays does not constitute a limitation to the approach presented here \cite{CGLW16,CKL20}, nor the approach is limited to the used of center-unstable manifolds formulas. 
In that respect,  more general parameterizations of the disregarded modes such as introduced in \cite{CLM19_closure} could be relevant to derive effective reduced GK systems able to handle more complex bifurcations \cite{CLM23}. 

Questions about the role of the annual cycle on this new route to complexity, in particular regarding its interactions with the underlying UPOs is a topic that should also lead to instructive dynamical insights in the tradition of the so-called slow-mode chaos approach  to ENSO variability \cite{chang1994interactions,jng94,TSCJ94,jin1996nino,Galanti_al00}.

Finally, we would like to mention that strong evidences about  the prominent role of UPOs in structuring atmospheric events, corresponding to zonal and blocking events in a low resolution quasi-geostrophic model, have been documented in \cite{lucarini2020new}. 
Our study pointed out the subtle interactions that UPOs may have with stochastic disturbances to produce a wealth of temporal events. We hope that more studies about the role of UPOs in the understanding of climate variability will be 
pursued in the future.

\section*{Acknowledgements}
The authors  thank Tetsushi Saburi for constructive comments on the GK method.   The authors also express their gratitude to the anonymous referees for their insightful comments that helped improve the presentation and discussion of the results. 
This work has been partially supported by the Office of Naval Research (ONR) Multidisciplinary University Research Initiative (MURI) grant N00014-20-1-2023, by the National Science Foundation grant DMS-2108856, and by the European Research Council (ERC) under the European Union's Horizon 2020 research and innovation program (grant agreement no. 810370).

\appendix

\section{Determination of $\mathcal{P}_N$ and $\mathcal{Q}_N$ in Eq.~\eqref{Eq_Galerkin_SS_6D}} \label{sect:coef_matrix_proof}
The calculations of $\mathcal{P}_N$ and $\mathcal{Q}_N$ are based on the RHS of \eqref{eq:A_SS}. $Q_N$ depends on the model parameters, while $\mathcal{P}_N$ is independent of those as explained below. 

First, the terms 
$$\frac{1}{\|\mathcal{K}_j\|_{\mathcal{E}}^2 } \sum_{n=0}^{N-1}  ( 1 - 3 T_{+}^2  - \alpha K_n(-1)),$$ when collected for each $j$ and $n$ allow us to build up $Q_N$ readily once  one recalls \eqref{eq:Kn_norm} and that due to \eqref{eq:Pn},
\be\label{Eq_K-1}
K_n(-1)=( n^2 + n + 1) (-1)^n,
\ee
since $L_n(-1)=(-1)^n$ from the properties of the Legendre polynomials.

The determination of $\mathcal{P}_N$ which results from the collection of the terms $ \sum_{k=0}^{n-1} a_{n,k} \left( \delta_{j,k} \|\mathcal{K}_j\|^2_{\mathcal{E}} - 1 \right )$, 
requires  the knowledge of the coefficients $a_{n,k}$. 
These coefficients arise in the expression of the derivatives of the Koornwinder polynomials in terms of these polynomials themselves; see \cite[Appendix~B]{CGLW16} for more details.  

The coefficients $a_{n,k}$ are then obtained by simply solving the triangular system \eqref{eq:algebraic} described in the following proposition.

\begin{prop_app} \label{prop:dPn}
The Koornwinder polynomial $K_n$ of degree $n$ defined in \eqref{eq:Pn} satisfies the differential relation
\be \label{eq:dPn}
\frac{\d K_n}{\d s}(s) = \sum_{k = 0}^{n-1} a_{n,k} K_k(s), \quad s \in (-1,1),
\ee
where the vector $\boldsymbol{a}_n$ made of the  $a_{n,k}$-coefficients,  
$$\boldsymbol{a}_n=(a_{n,0}, \cdots, a_{n,n-1})^T,$$
 solves the  following upper triangular system
\be \label{eq:algebraic}
\mathbf{T}\boldsymbol{a}_n = \boldsymbol{b}_n,
\ee
with $\mathbf{T}=(\mathbf{T}_{i,j})_{n\times n}$ and $\boldsymbol{b}_{n}=(b_{n,0}, \cdots, b_{n,n-1})^T$ given by
\bea \label{eq:algebraic_def}
\mathbf{T}_{i,j} & = \begin{cases}
0, & \; \text{ if } j < i,\\
i^2 + 1, & \; \text{ if } j = i,\\
-(2i+1), & \; \text{ if } j > i,
\end{cases} \; \quad \text{ where } \quad 0 \le i, j \le n-1, \\
b_{n,i} & = \begin{cases}
 -\frac{1}{2}(2i+1)(n+i+1)(n-i), & \text{ if $n+i$ is even}, \vspace{1em}\\
(n^2 + n)(2i+1) - \frac{i}{2}(n+i)(n-i+1) & \\
\hspace{2em} -\frac{1}{2}(i+1)(n-i-1)(n+i+2), & \text{ if $n+i$ is odd}.
\end{cases}
\eea

Finally, the rescaled Koornwinder polynomials  satisfy 
\be \label{eq:dKn_tau}
\frac{\d K^\tau_n}{\d \theta}(\theta) = \frac{2}{\tau} \sum_{k = 0}^{n-1} a_{n,k} K^\tau_k(\theta),  \quad \theta\in (-\tau,0).
\ee
\end{prop_app}

\begin{rem_app}
We note that the formula for $b_{n,i}$ that appeared in \cite{CGLW16} and in \cite{CKL20} contains two typos (a sign error and a factor 1/2 missing) that are here rectified in the formula \eqref{eq:algebraic_def}. We would like to express our gratitude to Tetsushi Saburi for pointing out these typos to us. We emphasize though that it is only a typo as the correct formulas were implemented in our codes, thus not questioning the diverse numerical results published in \cite{CGLW16}, \cite{CKL18_DDE}, \cite{CKL20}, and \cite{Chekroun_al22SciAdv}.
\end{rem_app}

\newcommand{\etalchar}[1]{$^{#1}$}
\providecommand{\bysame}{\leavevmode\hbox to3em{\hrulefill}\thinspace}
\providecommand{\MR}{\relax\ifhmode\unskip\space\fi MR }
\providecommand{\MRhref}[2]{%
  \href{http://www.ams.org/mathscinet-getitem?mr=#1}{#2}
}
\providecommand{\href}[2]{#2}


\begin{thebibliography}{ZGMPK03}

\bibitem[ABL{\etalchar{+}}22]{ando202215}
A.~And{\`o}, D.~Breda, D.~Liessi, S.~Maset, F.~Scarabel, and R.~Vermiglio,
  \emph{15 years or so of pseudospectral collocation methods for stability and
  bifurcation of delay equations}, Accounting for Constraints in Delay Systems,
  Springer, 2022, pp.~127--149.

\bibitem[ACN16]{Ashwin2016}
P.~Ashwin, S.~Coombes, and R.~Nicks, \emph{Mathematical frameworks for
  oscillatory network dynamics in neuroscience}, Journal of Mathematical
  Neuroscience \textbf{6} (2016), no.~2, 1--92.

\bibitem[AR23]{anikushin2023hidden}
M.~Anikushin and A.~Romanov, \emph{{Hidden and unstable periodic orbits as a
  result of homoclinic bifurcations in the Suarez--Schopf delayed oscillator
  and the irregularity of ENSO}}, Physica D: Nonlinear Phenomena \textbf{445}
  (2023), 133653.

\bibitem[Aut]{Auto}
\emph{{AUTO}}, \url{https://github.com/auto-07p/}, Accessed: 2023-09-09.

\bibitem[AWVC12]{ashwin2012tipping}
P.~Ashwin, S.~Wieczorek, R.~Vitolo, and P.~Cox, \emph{{Tipping points in open
  systems: Bifurcation, noise-induced and rate-dependent examples in the
  climate system}}, Phil. Trans. Roy. Soc. A \textbf{370} (2012), no.~1962,
  1166--1184.

\bibitem[BB78]{Banks_al78}
H.~T. Banks and J.~A. Burns, \emph{{Hereditary control problems: Numerical
  methods based on averaging approximations}}, SIAM J. Control Optim.
  \textbf{16} (1978), 169--208.

\bibitem[BCL{\etalchar{+}}17]{boers2017inverse}
N.~Boers, M.~D. Chekroun, H.~Liu, D.~Kondrashov, D.-D. Rousseau, A.~Svensson,
  M.~Bigler, and M.~Ghil, \emph{{Inverse stochastic--dynamic models for
  high-resolution Greenland ice core records}}, Earth Syst. Dynam. \textbf{8}
  (2017), 1171--1190.

\bibitem[BDG{\etalchar{+}}16]{breda2016pseudospectral}
D.~Breda, O.~Diekmann, M.~Gyllenberg, F.~Scarabel, and R.~Vermiglio,
  \emph{Pseudospectral discretization of nonlinear delay equations: new
  prospects for numerical bifurcation analysis}, SIAM Journal on Applied
  Dynamical Systems \textbf{15} (2016), 1--23.

\bibitem[BH89]{BH89}
D.~S. Battisti and A.~C. Hirst, \emph{{Interannual variability in a tropical
  atmosphere-ocean model: Influence of the basic state, ocean geometry and
  nonlinearity}}, J. Atmos. Sci. \textbf{46} (1989), 1687--1712.

\bibitem[BH05]{batstone2005characteristics}
C.~Batstone and H.H. Hendon, \emph{{Characteristics of stochastic variability
  associated with ENSO and the role of the MJO}}, Journal of Climate
  \textbf{18} (2005), no.~11, 1773--1789.

\bibitem[BK79]{banks1979spline}
H.~T. Banks and F.~Kappel, \emph{Spline approximations for functional
  differential equations}, Journal of Differential Equations \textbf{34}
  (1979), no.~3, 496--522.

\bibitem[BK98]{BK98}
W.-J. Beyn and W.~Kle\ss, \emph{{Numerical Taylor expansions of invariant
  manifolds in large dynamical systems}}, Numer. Math. \textbf{80} (1998),
  no.~1, 1--38.

\bibitem[BMV05]{breda2005pseudospectral}
D.~Breda, S.~Maset, and R.~Vermiglio, \emph{Pseudospectral differencing methods
  for characteristic roots of delay differential equations}, SIAM Journal on
  Scientific Computing \textbf{27} (2005), 482--495.

\bibitem[BNG97]{blanke1997estimating}
B.~Blanke, J.D. Neelin, and D.~Gutzler, \emph{{Estimating the effect of
  stochastic wind stress forcing on ENSO irregularity}}, Journal of Climate
  \textbf{10} (1997), no.~7, 1473--1486.

\bibitem[BRI84]{Banks_al84}
H.~T. Banks, I.~G. Rosen, and K.~Ito, \emph{{A spline based technique for
  computing Riccati operators and feedback controls in regulator problems for
  delay equations}}, SIAM J. Sci. Stat. Comput. \textbf{5} (1984), 830--855.

\bibitem[BTR07]{boutle2007nino}
I.~Boutle, R.H.S. Taylor, and R.A. R{\"o}mer, \emph{{El Ni{\~n}o and the
  delayed action oscillator}}, American Journal of Physics \textbf{75} (2007),
  15--24.

\bibitem[BVV14]{breda2014approximating}
D.~Breda and E.~Van~Vleck, \emph{{Approximating Lyapunov exponents and
  Sacker--Sell spectrum for retarded functional differential equations}},
  Numerische Mathematik \textbf{126} (2014), no.~2, 225--257.

\bibitem[BZ13]{bellenbook}
A.~Bellen and M.~Zennaro, \emph{{Numerical Methods for Delay Differential
  Equations}}, Numerical Mathematics and Scientific Computation, Oxford
  University Press, 2013.

\bibitem[CCH{\etalchar{+}}16]{chen2016diversity}
C.~Chen, M.A. Cane, N.~Henderson, D.~E. Lee, D.~Chapman, D.~Kondrashov, and
  M.~D. Chekroun, \emph{{Diversity, nonlinearity, seasonality, and memory
  effect in ENSO simulation and prediction using empirical model reduction}},
  Journal of Climate \textbf{29} (2016), no.~5, 1809--1830.

\bibitem[CCHT19]{cao2019mathematical}
Y.~Cao, M.~D. Chekroun, A.~Huang, and R.~Temam, \emph{{Mathematical Analysis of
  the Jin-Neelin Model of El Ni\~no-Southern-Oscillation}}, Chinese Annals of
  Mathematics, Series B \textbf{40} (2019), no.~1, 1--38.

\bibitem[CFY22]{chen2022multiscale}
N.~Chen, X.~Fang, and J.-Y. Yu, \emph{{A multiscale model for El Ni{\~n}o
  complexity}}, npj Climate and Atmospheric Science \textbf{5} (2022), 16.

\bibitem[CGLW16]{CGLW16}
M.~D. Chekroun, M.~Ghil, H.~Liu, and S.~Wang, \emph{{Low-dimensional Galerkin
  approximations of nonlinear delay differential equations}}, Disc. Cont. Dyn.
  Sys. A \textbf{36} (2016), no.~8, 4133--4177.

\bibitem[CGN18]{CGN18}
M.D. Chekroun, M.~Ghil, and J.~D. Neelin, \emph{{Pullback attractor crisis in a
  delay differential ENSO model}}, {Advances in Nonlinear Geosciences}
  (A.~Tsonis, ed.), Springer, 2018, pp.~1--33.

\bibitem[CK17]{chekroun2017data}
M.~D. Chekroun and D.~Kondrashov, \emph{{Data-adaptive harmonic spectra and
  multilayer Stuart-Landau models}}, Chaos \textbf{27} (2017), no.~9, 093110.

\bibitem[CKG11]{CKG11}
M.~D. Chekroun, D.~Kondrashov, and M.~Ghil, \emph{{Predicting stochastic
  systems by noise sampling, and application to the El Ni\~no-Southern
  Oscillation}}, {Proc. Natl. Acad. Sci.} \textbf{108} (2011), 11766--11771.

\bibitem[CKL18]{CKL18_DDE}
M.~D. Chekroun, A.~Kr\"oner, and H.~Liu, \emph{Galerkin approximations for the
  optimal control of nonlinear delay differential equations},
  {Hamilton-Jacobi-Bellman Equations. Numerical Methods and Applications in
  Optimal Control. D. Kalise, K. Kunisch, and Z. Rao (Eds.)}, vol.~21, Berlin,
  Boston: De Gruyter, 2018, pp.~61--96.

\bibitem[CKL20]{CKL20}
M.~D. Chekroun, I.~Koren, and H.~Liu, \emph{{Efficient reduction for diagnosing
  Hopf bifurcation in delay differential systems: Applications to cloud-rain
  models}}, Chaos \textbf{40} (2020), no.~8, 053130,
  \href{https://doi.org/10.1063/5.0004697}{doi:10.1063/5.0004697}.

\bibitem[CKLL22]{Chekroun_al22SciAdv}
M.~D. Chekroun, I.~Koren, H.~Liu, and H.~Liu, \emph{Generic generation of
  noise-driven chaos in stochastic time delay systems: Bridging the gap with
  high-end simulations}, Science Advances \textbf{8} (2022), no.~46, eabq7137,
  \href{https://doi.org/10.1126/sciadv.abq7137}{doi.org/10.1126/sciadv.abq7137}.

\bibitem[CLM20]{CLM19_closure}
M.~D. Chekroun, H.~Liu, and J.~C. McWilliams, \emph{{Variational approach to
  closure of nonlinear dynamical systems: Autonomous case}}, Journal of
  Statistical Physics \textbf{179} (2020), 1073--1160.

\bibitem[CLM23]{CLM23}
M.D. Chekroun, H.~Liu, and J.C. McWilliams, \emph{Optimal parameterizing
  manifolds for anticipating tipping points and higher-order critical
  transitions}, Chaos \textbf{33} (2023), 093126.

\bibitem[CLMW23]{chekroun2023transitions}
M.~D. Chekroun, H.~Liu, J.~C. McWilliams, and S.~Wang, \emph{{Transitions in
  stochastic non-equilibrium systems: Efficient reduction and analysis}},
  Journal of Differential Equations \textbf{346} (2023), 145--204,
  \href{https://doi.org/10.1016/j.jde.2022.11.025}{doi.org/10.1016/j.jde.2022.11.025}.

\bibitem[CLW15a]{CLW15_vol1}
M.~D. Chekroun, H.~Liu, and S.~Wang, \emph{{Approximation of Stochastic
  Invariant Manifolds: Stochastic Manifolds for Nonlinear SPDEs I}}, Springer
  Briefs in Mathematics, Springer, New York, 2015.

\bibitem[CLW15b]{CLW15_vol2}
\bysame, \emph{{Stochastic Parameterizing Manifolds and Non-Markovian Reduced
  Equations: Stochastic Manifolds for Nonlinear SPDEs II}}, Springer Briefs in
  Mathematics, Springer, New York, 2015.

\bibitem[CME18]{EU_Nino34}
CMEMS, \emph{{Ni\~no 3.4 Sea Surface Temperature time series from Reanalysis}},
  E.U.~Copernicus Marine Service Information. Marine Data Store (MDS) (2018),
  \href{https://data.marine.copernicus.eu/product/GLOBAL_OMI_CLIMVAR_enso_sst_area_averaged_anomalies/description}{doi:10.48670/moi-00219}.

\bibitem[CMT18]{chen2018observations}
N.~Chen, A.~J. Majda, and S.~Thual, \emph{{Observations and mechanisms of a
  simple stochastic dynamical model capturing El Ni{\~n}o diversity}}, Journal
  of Climate \textbf{31} (2018), 449--471.

\bibitem[Cra91]{Crawford91}
J.~D. Crawford, \emph{Introduction to bifurcation theory}, Reviews of Modern
  Physics \textbf{63} (1991), no.~4, 991--1036.

\bibitem[CWLJ94]{chang1994interactions}
P.~Chang, B.~Wang, T.~Li, and L.~Ji, \emph{{Interactions between the seasonal
  cycle and the Southern Oscillation-Frequency entrainment and chaos in a
  coupled ocean-atmosphere model}}, Geophysical Research Letters \textbf{21}
  (1994), no.~25, 2817--2820.

\bibitem[CZ95]{curtain1995}
R.F. Curtain and H.~Zwart, \emph{{An Introduction to Infinite-Dimensional
  Linear Systems Theory}}, vol.~21, Springer, 1995.

\bibitem[CZ23]{chen2023rigorous}
N.~Chen and Y.~Zhang, \emph{{Rigorous derivation of stochastic conceptual
  models for the El Ni{\~n}o-Southern Oscillation from a spatially-extended
  dynamical system}}, Physica D: Nonlinear Phenomena \textbf{453} (2023),
  133842.

\bibitem[DCP{\etalchar{+}}21]{dieppois2021enso}
B.~Dieppois, A.~Capotondi, B.~Pohl, K.P. Chun, P.-A. Monerie, and J.~Eden,
  \emph{{ENSO diversity shows robust decadal variations that must be captured
  for accurate future projections}}, Communications Earth \& Environment
  \textbf{2} (2021), no.~1, 212.

\bibitem[DDMG16]{detrixhe2016fast}
M.~Detrixhe, M.~Doubeck, J.~Moehlis, and F.~Gibou, \emph{{A fast Eulerian
  approach for computation of global isochrons in high dimensions}}, SIAM
  Journal on Applied Dynamical Systems \textbf{15} (2016), no.~3, 1501--1527.

\bibitem[DGK{\etalchar{+}}08]{dhooge2008new}
A.~Dhooge, W.~Govaerts, Yu.~A. Kuznetsov, H.~G.~E. Meijer, and B.~Sautois,
  \emph{{New features of the software MatCont for bifurcation analysis of
  dynamical systems}}, Mathematical and Computer Modelling of Dynamical Systems
  \textbf{14} (2008), 147--175.

\bibitem[EL97]{eckert1997predictability}
C.~Eckert and M.~Latif, \emph{{Predictability of a stochastically forced hybrid
  coupled model of El Ni{\~n}o}}, Journal of Climate \textbf{10} (1997), no.~7,
  1488--1504.

\bibitem[ELR02]{Engelborghs_al02}
K.~Engelborghs, T.~Luzyanina, and D.~Roose, \emph{{Numerical bifurcation
  analysis of delay differential equations using DDE-BIFTOOL}}, ACM
  Transactions on Mathematical Software \textbf{28} (2002), 1--21.

\bibitem[EvP04]{EvP04}
T.~Eirola and J.~von Pfaler, \emph{{Numerical Taylor expansions for invariant
  manifolds}}, Numer. Math. \textbf{99} (2004), no.~1, 25--46.

\bibitem[Fed02]{fedorov2002response}
A.V. Fedorov, \emph{The response of the coupled tropical ocean--atmosphere to
  westerly wind bursts}, Quarterly Journal of the Royal Meteorological Society
  \textbf{128} (2002), no.~579, 1--23.

\bibitem[FPS18]{feudel2018multistability}
U.~Feudel, A.~N. Pisarchik, and K.~Showalter, \emph{{Multistability and
  tipping: From mathematics and physics to climate and brain--Minireview and
  preface to the focus issue}}, Chaos \textbf{28} (2018), 033501.

\bibitem[FQS{\etalchar{+}}19]{falkena2019derivation}
S.~K.~J. Falkena, C.~Quinn, J.~Sieber, J.~Frank, and H.~A. Dijkstra,
  \emph{{Derivation of delay equation climate models using the Mori-Zwanzig
  formalism}}, Proceedings of the Royal Society A \textbf{475} (2019),
  no.~2227, 20190075.

\bibitem[Gar04]{gardinerbook}
C.~W. Gardiner, \emph{{Handbook of Stochastic Methods}}, 3rd ed.,
  Springer-Verlag Berlin, 2004.

\bibitem[GCS15]{GCStep15}
M.~Ghil, M.~D. Chekroun, and G.~Stepan, \emph{A collection on {`Climate
  Dynamics: Multiple Scales and Memory Effects', Editorial}}, R. Soc. Proc. A
  \textbf{471} (2015), 20150097.

\bibitem[GHJM98]{gammaitoni1998stochastic}
L.~Gammaitoni, P.~H{\"a}nggi, P.~Jung, and F.~Marchesoni, \emph{Stochastic
  resonance}, Reviews of Modern Physics \textbf{70} (1998), 223--287.

\bibitem[Gri08]{gritsun2008unstable}
A.S. Gritsun, \emph{Unstable periodic trajectories of a barotropic model of the
  atmosphere}, Russ. J. Numer. Anal. Math. Modelling \textbf{43} (2008),
  345--367.

\bibitem[Gri13]{gritsun2013statistical}
A.~Gritsun, \emph{Statistical characteristics, circulation regimes and unstable
  periodic orbits of a barotropic atmospheric model}, Philosophical
  Transactions of the Royal Society A: Mathematical, Physical and Engineering
  Sciences \textbf{371} (2013), no.~1991, 20120336.

\bibitem[GT00]{Galanti_al00}
E.~Galanti and E.~Tziperman, \emph{{ENSO's phase locking to the seasonal cycle
  in the fast-SST, fast-wave, and mixed-mode regimes}}, J. Atmos. Sci.
  \textbf{57} (2000), 2936--2950.

\bibitem[Guc75]{guckenheimer1975isochrons}
J.~Guckenheimer, \emph{Isochrons and phaseless sets}, Journal of Mathematical
  Biology \textbf{1} (1975), no.~3, 259--273.

\bibitem[GZT08]{ghil2008delay}
M.~Ghil, I.~Zaliapin, and S.~Thompson, \emph{{A delay differential model of
  ENSO variability: parametric instability and the distribution of extremes}},
  Nonlinear Processes in Geophysics \textbf{15} (2008), no.~3, 417--433.

\bibitem[Hal88]{Hale88}
J.~K. Hale, \emph{{A}symptotic {B}ehavior of {D}issipative {S}ystems},
  Mathematical Surveys and Monographs, vol.~25, American Mathematical Society,
  Providence, RI, 1988.

\bibitem[HCF{\etalchar{+}}16]{haro2016parameterization}
{\`A}.~Haro, M.~Canadell, J.-L. Figueras, A.~Luque, and J.-M. Mondelo,
  \emph{{The Parameterization Method for Invariant Manifolds:From Rigorous
  Results to Effective Computations}}, vol. 195, Springer-Verlag, Berlin, 2016.

\bibitem[Hen81]{Hen81}
D.~Henry, \emph{Geometric {T}heory of {S}emilinear {P}arabolic {E}quations},
  Lecture Notes in Mathematics, vol. 840, Springer-Verlag, Berlin, 1981.

\bibitem[HL84]{Horsthemke_book}
W.~Horsthemke and R.~Lefever, \emph{{Noise-induced Transitions}}, Springer
  Series in Synergetics, vol.~15, Springer-Verlag, Berlin, 1984, Theory and
  applications in physics, chemistry, and biology.

\bibitem[HVL93]{Hale_Lunel93}
J.~K. Hale and S.~M. Verduyn~Lunel, \emph{{Introduction to
  Functional-Differential Equations}}, Applied Mathematical Sciences, vol.~99,
  Springer-Verlag, New York, 1993.

\bibitem[IT86]{Ito_Teglas86}
K.~Ito and R.~Teglas, \emph{Legendre-tau approximations for
  functional-differential equations}, SIAM J. Control Optim. \textbf{24}
  (1986), no.~4, 737--759.

\bibitem[JN93]{Jin_al93_part1}
F.-F. Jin and J.~D. Neelin, \emph{{Modes of interannual tropical
  ocean-atmosphere interaction-A unified view. Part I: Numerical results}},
  Journal of the atmospheric sciences \textbf{50} (1993), no.~21, 3477--3503.

\bibitem[JNG94]{jng94}
F.-F. Jin, J.~D. Neelin, and M.~Ghil, \emph{{El Ni\~no on the Devil's
  Staircase: Annual subharmonic steps to chaos}}, Science \textbf{274} (1994),
  70--72.

\bibitem[JNG96]{jin1996nino}
F.-F. Jin, J.D. Neelin, and M.~Ghil, \emph{{El Ni{\~n}o/Southern Oscillation
  and the annual cycle: Subharmonic frequency-locking and aperiodicity}},
  Physica D: Nonlinear Phenomena \textbf{98} (1996), no.~2-4, 442--465.

\bibitem[Kap86]{Kappel86}
F.~Kappel, \emph{Semigroups and delay equations}, {Semigroups, Theory and
  Applications, Vol.~II (Trieste, 1984)}, Pitman Res. Notes Math. Ser., vol.
  152, Longman Sci. Tech., Harlow, 1986, pp.~136--176.

\bibitem[KC18]{kondrashov2018data}
D.~Kondrashov and M.~D. Chekroun, \emph{Data-adaptive harmonic analysis and
  modeling of solar wind-magnetosphere coupling}, Journal of Atmospheric and
  Solar-Terrestrial Physics \textbf{177} (2018), 179--189.

\bibitem[KCB18]{KCB18}
D.~Kondrashov, M.~D. Chekroun, and P.~Berloff, \emph{Multiscale
  {S}tuart-{L}andau emulators: {A}pplication to wind-driven ocean gyres},
  Fluids \textbf{3} (2018), no.~1, 21.

\bibitem[KCG18]{MASIE_paper}
D.~Kondrashov, M.D. Chekroun, and M.~Ghil, \emph{Data-adaptive harmonic
  decomposition and prediction of {A}rctic sea ice extent}, Dynamics and
  Statistics of the Climate System \textbf{3} (2018), no.~1, 1--23.

\bibitem[KCYG18]{KCYG17}
D.~Kondrashov, M.~D. Chekroun, X.~Yuan, and M.~Ghil, \emph{{Data-adaptive
  harmonic decomposition and stochastic modeling of Arctic sea ice}}, Advances
  in Nonlinear Geosciences ({A. Tsonis}, ed.), Springer, 2018, pp.~179--205.

\bibitem[KF11]{koren2011aerosol}
I.~Koren and G.~Feingold, \emph{Aerosol--cloud--precipitation system as a
  predator-prey problem}, Proc. Natl. Acad. Sci. USA \textbf{108} (2011),
  no.~30, 12227--12232.

\bibitem[KKD19]{keane2019effect}
A.~Keane, B.~Krauskopf, and H.~A. Dijkstra, \emph{{The effect of state
  dependence in a delay differential equation model for the El Ni{\~n}o
  Southern Oscillation}}, Philosophical Transactions of the Royal Society A
  \textbf{377} (2019), no.~2153, 20180121.

\bibitem[KKP15]{keane2015delayed}
A.~Keane, B.~Krauskopf, and C.~Postlethwaite, \emph{{Delayed feedback versus
  seasonal forcing: Resonance phenomena in an El Ni{\~n}o Southern Oscillation
  model}}, SIAM Journal on Applied Dynamical Systems \textbf{14} (2015), no.~3,
  1229--1257.

\bibitem[KKP16]{keane2017}
A.~Keane, B.~Krauskopf, and C.~M. Postlethwaite, \emph{{Investigating irregular
  behavior in a model for the El Ni{\~n}o Southern Oscillation with positive
  and negative delayed feedback}}, SIAM J. Appl. Dyn. Syst. \textbf{15} (2016),
  1656--1689.

\bibitem[KKP17]{keane2017climate}
\bysame, \emph{Climate models with delay differential equations}, Chaos: An
  Interdisciplinary Journal of Nonlinear Science \textbf{27} (2017), no.~11,
  114309.

\bibitem[KKRG05]{kondrashov2005hierarchy}
D.~Kondrashov, S.~Kravtsov, A.~W. Robertson, and M.~Ghil, \emph{{A hierarchy of
  data-based ENSO models}}, Journal of Climate \textbf{18} (2005), 4425--4444.

\bibitem[Knu]{Knut}
\emph{{KNUT}}, \url{http://rs1909.github.io/knut/}, Accessed: 2023-09-09.

\bibitem[Koo84]{Koo84}
T.~H. Koornwinder, \emph{Orthogonal polynomials with weight function
  $(1-x)^\alpha(1+ x) ^\beta + {M} \delta (x+ 1)+ {N} \delta (x-1)$}, Canad.
  Math. Bull. \textbf{27} (1984), 205--214.

\bibitem[KS78]{kappel1978autonomous}
F.~Kappel and W.~Schappacher, \emph{Autonomous nonlinear functional
  differential equations and averaging approximations}, Nonlinear Analysis:
  Theory, Methods \& Applications \textbf{2} (1978), no.~4, 391--422.

\bibitem[KS14]{krauskopf2014bifurcation}
B.~Krauskopf and J.~Sieber, \emph{{Bifurcation analysis of delay-induced
  resonances of the El-Ni{\~n}o Southern Oscillation}}, Proceedings of the
  Royal Society A: Mathematical, Physical and Engineering Sciences \textbf{470}
  (2014), no.~2169, 20140348.

\bibitem[KTF17]{koren2017exploring}
I.~Koren, E.~Tziperman, and G.~Feingold, \emph{Exploring the nonlinear cloud
  and rain equation}, Chaos: An Interdisciplinary Journal of Nonlinear Science
  \textbf{27} (2017), no.~1, 013107.

\bibitem[Kue11]{kuehn2011mathematical}
Ch. Kuehn, \emph{{A mathematical framework for critical transitions:
  Bifurcations, fast--slow systems and stochastic dynamics}}, Physica D
  \textbf{240} (2011), no.~12, 1020--1035.

\bibitem[LG20]{lucarini2020new}
V.~Lucarini and A.~Gritsun, \emph{A new mathematical framework for atmospheric
  blocking events}, Climate Dynamics \textbf{54} (2020), no.~1-2, 575--598.

\bibitem[LTW{\etalchar{+}}17]{BAMS_2017observing}
M.~L. L'Heureux, K.~Takahashi, A.~B. Watkins, A.~G. Barnston, E.~J. Becker,
  T.~E. Di~Liberto, F.~Gamble, J.~Gottschalck, M.~S. Halpert, B.~Huang, et~al.,
  \emph{{Observing and predicting the 2015/16 El Ni{\~n}o}}, Bulletin of the
  American Meteorological Society \textbf{98} (2017), no.~7, 1363--1382.

\bibitem[Mat]{Matcont}
\emph{{MatCont}}, \url{https://sourceforge.net/projects/matcont/}, Accessed:
  2023-09-09.

\bibitem[MRMM14]{mauroy2014global}
A.~Mauroy, B.~Rhoads, J.~Moehlis, and I.~Mezic, \emph{Global isochrons and
  phase sensitivity of bursting neurons}, SIAM Journal on Applied Dynamical
  Systems \textbf{13} (2014), no.~1, 306--338.

\bibitem[MW05]{MW05}
T.~Ma and S.~Wang, \emph{Bifurcation {T}heory and {A}pplications}, World
  Scientific Series on Nonlinear Science. Series A: Monographs and Treatises,
  vol.~53, World Scientific Publishing Co. Pte. Ltd., Hackensack, NJ, 2005.

\bibitem[MW14]{MW14}
\bysame, \emph{{Phase Transition Dynamics}}, Springer, 2014.

\bibitem[NBH{\etalchar{+}}98]{Neelin_al98}
J.~D. Neelin, D.~S. Battisti, A.~C. Hirst, F.-F. Jin, Y.~Wakata, T.~Yamagata,
  and S.~E. Zebiak, \emph{{ENSO theory}}, J. Geophys. Res. Oceans \textbf{103}
  (1998), 14261--14290.

\bibitem[NJ93]{Jin_al93_part2}
J.~D. Neelin and F.-F. Jin, \emph{{Modes of interannual tropical
  ocean-atmosphere interaction-a unified view. Part II: Analytical results in
  the weak-coupling limit}}, Journal of the atmospheric sciences \textbf{50}
  (1993), no.~21, 3504--3522.

\bibitem[PS95]{penland1995optimal}
C.~Penland and P.~D. Sardeshmukh, \emph{The optimal growth of tropical sea
  surface temperature anomalies}, Journal of climate \textbf{8} (1995), no.~8,
  1999--2024.

\bibitem[RCC{\etalchar{+}}14]{roques2014parameter}
L.~Roques, M.~D. Chekroun, M.~Cristofol, S.~Soubeyrand, and M.~Ghil,
  \emph{Parameter estimation for energy balance models with memory},
  Proceedings of the Royal Society A: Mathematical, Physical and Engineering
  Sciences \textbf{470} (2014), no.~2169, 20140349.

\bibitem[RN00]{roulston2000response}
M.S. Roulston and J.D. Neelin, \emph{{The response of an ENSO model to climate
  noise, weather noise and intraseasonal forcing}}, Geophysical Research
  Letters \textbf{27} (2000), no.~22, 3723--3726.

\bibitem[SBB{\etalchar{+}}09]{scheffer2009early}
M.~Scheffer, J.~Bascompte, W.~A. Brock, V.~Brovkin, S.~R. Carpenter, V.~Dakos,
  H.~Held, E.~H. Van~Nes, M.~Rietkerk, and G.~Sugihara, \emph{Early-warning
  signals for critical transitions}, Nature \textbf{461} (2009), 53--59.

\bibitem[Sel85]{sell1985smooth}
G.~R. Sell, \emph{Smooth linearization near a fixed point}, American Journal of
  Mathematics (1985), 1035--1091.

\bibitem[SEL{\etalchar{+}}14]{sieber7144dde}
J.~Sieber, K.~Engelborghs, T.~Luzyanina, G.~Samaey, and D.~Roose,
  \emph{{DDE-BIFTOOL v. 3.1.1 Manual---Bifurcation analysis of delay
  differential equations}}, arXiv preprint arXiv:1406.7144 (2014).

\bibitem[SS88]{Suarez_al88}
M.~J. Suarez and P.~S. Schopf, \emph{{A delayed action oscillator for ENSO}},
  J. Atmos. Sci. \textbf{45} (1988), 3283--3287.

\bibitem[SY02]{SY02}
G.~Sell and Y.~You, \emph{Dynamics of {E}volutionary {E}quations}, Applied
  Mathematical Sciences, vol. 143, Springer-Verlag, New York, 2002.

\bibitem[TCDN20]{Tantet_al_ENSO}
A.~Tantet, M.D. Chekroun, H.~Dijkstra, and J.~D. Neelin,
  \emph{{Ruelle-Pollicott Resonances of Stochastic Systems in Reduced State
  Space. Part III: Application to the Cane-Zebiak Model of the El Nino-Southern
  Oscillation}}, J.~Stat.~Phys. \textbf{179} (2020), 1449--1474,
  \href{https://doi.org/10.1007/s10955-019-02444-8}{doi:
  10.1007/s10955-019-02444-8}.

\bibitem[TMCS16]{thual2016simple}
S.~Thual, A.~J. Majda, N.~Chen, and S.~N. Stechmann, \emph{{Simple stochastic
  model for El Ni{\~n}o with westerly wind bursts}}, Proc. Natl. Acad. Sci. USA
  \textbf{113} (2016), 10245--10250.

\bibitem[TSCJ94]{TSCJ94}
E.~Tziperman, L.~Stone, M.~A. Cane, and H.~Jarosh, \emph{{{El Ni{\~n}o chaos:
  Overlapping of resonances between the seasonal cycle and the Pacific
  ocean-atmosphere oscillator}}}, Science \textbf{264} (1994), 72--74.

\bibitem[Vya12]{Vyasarayani12}
C.~P. Vyasarayani, \emph{Galerkin approximations for higher order delay
  differential equations}, J. Comput. Nonlinear Dynamics \textbf{7} (2012),
  031004.

\bibitem[WC05]{Wahi_al05}
P.~Wahi and A.~Chatterjee, \emph{Galerkin projections for delay differential
  equations}, ASME. J. Dyn. Syst., Meas., Control \textbf{127} (2005), 80--87.

\bibitem[ZGMPK03]{zavala2003response}
J.~Zavala-Garay, A.M. Moore, C.L. Perez, and R.~Kleeman, \emph{{The response of
  a coupled model of ENSO to observed estimates of stochastic forcing}},
  Journal of Climate \textbf{16} (2003), no.~17, 2827--2842.

\end{thebibliography}
\end{document}